\documentclass[review]{elsarticle}
\usepackage[a4paper, margin=1in]{geometry}
\usepackage{amsmath}
\usepackage{stackengine}
\usepackage{bigints}
\stackMath
\usepackage{lineno}
\modulolinenumbers[5]

\usepackage{booktabs}
\usepackage{mathptmx}      
\usepackage{latexsym}
\usepackage{graphicx}
\usepackage{amssymb}
\usepackage{amsmath,bm}
\usepackage{times}
\usepackage{lipsum}
\usepackage{caption}
\usepackage{subcaption}
\usepackage{rotating}
\usepackage{xcolor}
\usepackage{mathtools}
\usepackage{tikz}
\usepackage[utf8]{inputenc}
\usetikzlibrary{arrows.meta, shapes.geometric}
\usepackage{hyperref}
\hypersetup{
    colorlinks=true,
    linkcolor=black,
    filecolor=black,      
    urlcolor=black,
    citecolor=black,
    linkbordercolor=black,
}

\usepackage{titlesec}
\titleformat{\section}
{\normalfont\Large\bfseries}{\thesection}{1em}{}
\titleformat{\subsection}
{\normalfont\large\bfseries}{\thesubsection}{1em}{}

\usepackage{amsmath}
\usepackage{array}
\usepackage{bm}
\usepackage{longtable}

\journal{Applied Numerical Mathematics}

\makeatletter
\def\ps@pprintTitle{}
\makeatother

\begin{document}

\begin{frontmatter}

\journal{Applied Numerical Mathematics}

\title{{Finite Element Simulation of Modified \\Poisson–Nernst–Planck/Navier–Stokes Model for \\ Compressible Electrolytes under Mechanical Equilibrium}}

\tnotetext[1]{Accepted for publication in \textit{Applied Numerical Mathematics}.}


\author[1,2]{Ankur Ankur}
\ead{Ankur@ma.iitr.ac.in}
\author[2]{Ram Jiwari}
\ead{ram.jiwari@ma.iitr.ac.in}
\author[3]{Satyvir Singh \corref{cor1}}
\ead{singh@acom.rwth-aachen.de}
\cortext[cor1]{Corresponding author}
\address[1]{Mathematics Area, SISSA, International School for Advanced Studies, via Bonomea 265, Trieste 34136, Italy}
\address[2]{Department of Mathematics, Indian Institute of Technology Roorkee, Roorkee 247667, India}
\address[3]{Applied and Computational Mathematics, RWTH Aachen University, Aachen 52062, Germany}


\begin{abstract} 
This work presents a finite element method for a modified Poisson–Nernst–Planck/Navier–Stokes (PNP/NS) model under the mechanical equilibrium, developed for compressible electrolytes. The modification is based on the new model proposed by Dreyer, Guhlke and M$\ddot{u}$ller \cite{dreyer2013overcoming}, where the diffusion flux in the classical PNP system is replaced with an implicitly involved new diffusion flux, leading to fractional nonlinearity. He and Sun \cite{Mingyan} previously developed a numerical approach for another type of modification, where the Poisson equation in the PNP system was substituted with a fourth-order elliptic equation. 
Another key contribution of this work is the reduction of the equilibrium system to a modified Poisson–Boltzmann system.
The proposed numerical scheme is capable of handling both compressible and incompressible regimes by employing a bulk modulus parameter, which governs the fluid's compressibility and enables seamless transition between these regimes. To emphasize practical relevance, we discuss the implications of compressible electrolytes in the context of double-layer capacitance behavior.
We also conduct numerical simulations over various domains to demonstrate its applicability under various operating conditions, including temperature fluctuations and variations in the bulk modulus. The numerical results validate the accuracy and robustness of our computational scheme and demonstrate that the observed limiting behavior for the incompressible regime aligns with the theoretical trends anticipated by  Dreyer et al. \cite{dreyer2013overcoming}.
\end{abstract}


\end{frontmatter}

\section{\Large Introduction}
\label{Sec:1}
\subsection{\large \textbf{Scope}}
The popular Poisson–Nernst–Planck/Navier–Stokes (PNP/NS) model is derived from the classical theory of electrokinetics. It is utilized to mathematically describe the dynamical properties of electrically charged fluids, the motions of ionized particles  and their interactions with the surrounding fluid and electric fields  \cite{jerome}. Ion concentrations are governed by convection–diffusion–reaction equations, commonly known as Nernst–Planck (NP) equations, the diffusive nature of the electrostatic potential is described by the Poisson equation, while the Navier-Stokes (NS) equations describe the dynamics of the fluids, neglecting the gravitational and magnetic forces. The PNP/NS coupling thus efficiently describes electro-chemical and fluid-mechanical transport across cellular environments, rendering this model applicable to a broad spectrum of phenomena, including drug delivery into biomembranes, electrokinetic flows in electrophysiology, semiconductors and many others  (see, e.g., \cite{hu2005, Cioffi2006, Wang2017, Lu2010, Jerome2008} and
the references therein). A visual representation of simple electrochemical cell with zinc and copper electrodes in NaCl electrolyte solution is depicted in Figure \ref{Fig:1}.

The limitations of the classical electrokinetic theory have long been recognized, leading to the development of new formulations \cite{ Kornyshev1981Conductivity, Tresset2008generalized, Landstorfer2011, sparnaay1958corrections,Mansoori1971equilibrium, dreyer2013overcoming}. Notably, the model proposed by Nernst and Planck \cite{Nernst1888, Nernst1889, Planck1890one, Planck1890two} neglects the mass conservation of the entire mixture and lacks coupling to the momentum balance equation. Subsequent models pertaining to this issue were similarly inconsistent with the theory of non-equilibrium thermodynamics \cite{dreyer2013overcoming}, which emerged more than eight decades ago.
As a result, they often predicted excessive pressure in the boundary layers, necessitating the extension of the classical electrokinetic theory.

Recently, Dreyer, Guhlke and M$\ddot{u}$ller \cite{dreyer2013overcoming} addressed these inherent deficiencies, leading to a new set of equations to describe the ionic concentrations for a mixture with $N-1$ ionic species and a neutral solvent: 
\begin{equation}
	\frac{\partial (m_{\alpha}n_{\alpha})}{\partial t} + \nabla\cdot(m_{\alpha}n_{\alpha}\bm{v}+\bm{J}_{\alpha})=0, \quad \alpha\;\in\;\{1,\ldots,N-1\}, \label{Eq19nn}
\end{equation} 
having diffusion fluxes 
\begin{equation}
\label{Eq8nn}
\bm{J}_{\alpha}= - \sum_{\beta=1}^{N-1}H_{\alpha\beta} \bigg(\frac{\nabla\mu_{\beta}-\nabla\mu_N}{T}+ \frac{e_0}{T}\bigg(\frac{z_{\beta}}{m_{\beta}}-\frac{z_{N}}{m_{N}}\bigg)\nabla\varphi\bigg), \quad\alpha \in \{1, 2, \ldots N-1\}, 
\end{equation}
with
\begin{equation}
\label{Eq9nn}
\mu_{\alpha}= g_{\alpha}^R+\frac{kT}{m_{\alpha}}\mbox{ln}\bigg(\frac{n_\alpha}{n}\bigg)+\frac{K}{m_{\alpha}n^R}\mbox{ln}\bigg(1+\frac{p-p^R}{K}\bigg) , \quad\quad\alpha \in \{1, 2, \ldots N\}. 
\end{equation}
\begin{equation}
\label{Eq10nn}
p= p^R+K\bigg(\frac{n}{n^R}-1\bigg), 
\end{equation}
 where,  $n_{\alpha}$ ($\alpha=1, 2, \ldots N$) are the number densities of ionic species with mass $m_{\alpha}$ and charge $z_{\alpha}$; $n$ is the total number density; $p, \varphi$, $\mu_{\alpha}, \bm{v}$, $T$, $g_{\alpha}$, $k$, $K$,$H_{\alpha\beta}$ and $e_0$ are pressure, electrostatic potential, chemical potential, velocity, temperature, specific Gibbs energies, Boltzmann constant, bulk modulus, positive definite kinetic matrix and elementary charge (constant), respectively and $R$ indicates the reference state. Although the complete mathematical modeling of the modified PNP/NS system has been addressed in \cite{dreyer2013overcoming}, it is reproduced in \ref{Sec:2app} for convenience, while a table summarizing all the notations used in this article can be found in \ref{Notations}.
 
 \vspace{0.2cm}
 
In contrast, the diffusion flux in the classical Nernst and Planck model is given by (e.g., see equation $20$ in \cite{dreyer2013overcoming})
\begin{equation}
\label{Eq8nnn}
\bm{J}_{\alpha}= - M_{\alpha}^{NP}(kT\nabla n_{\alpha}+z_{\alpha}e_0n_{\alpha}\nabla\varphi), \quad\alpha \in \{1, 2, \ldots N\}, 
\end{equation}
This form of the diffusion flux has been adopted in many other models, with 
$\alpha$ typically ranging up to $\{1, 2, \ldots N-1\}$. Fuhrmann \cite{Fuhrmann2015} provides a comparative analysis of several such modified formulations for the solutions of the PNP/NS system for incompressible electrolytes.


\begin{figure}[h!]
\centering

\begin{tikzpicture}[scale=0.8,
    cation/.style={circle,fill=blue!70,minimum size=5mm,inner sep=0pt},
    anion/.style={circle,fill=red!70,minimum size=5mm,inner sep=0pt},
    zincion/.style={circle,fill=green!70,minimum size=5mm,inner sep=0pt}]

\draw[thick, fill=blue!5] (0,0) rectangle (8,4);
\node at (4,4.3) {\small NaCl Electrolyte Solution};

\fill[gray!60] (0.5,0.5) rectangle (1.0,3.5);
\node[rotate=90] at (0.75,1.9) {Zn};
\node at (0.75,3.8) {$-$};

\fill[brown!60] (7.0,0.5) rectangle (7.5,3.5);
\node[rotate=90] at (7.25,1.9) {Cu};
\node at (7.25,3.8) {$+$};

\node[cation] at (2.5,1.2) {\tiny Na$^+$};
\node[cation] at (3.0,2.0) {\tiny Na$^+$};
\node[cation] at (4.0,0.8) {\tiny Na$^+$};
\node[cation] at (5.5,3.2) {\tiny Na$^+$};
\node[cation] at (3.8,1.7) {\tiny Na$^+$};
\node[cation] at (4.8,1.7) {\tiny Na$^+$};
\node[cation] at (6.4,1.7) {\tiny Na$^+$};
\draw[->,blue] (2.85,1.2) -- (3.4,1.2) node[midway,above] {\tiny Na$^+$};
\draw[->,blue] (5.85,3.2) -- (6.4,3.2) node[midway,above] {\tiny Na$^+$};

\node[anion] at (6.0,2.5) {\tiny Cl$^-$};
\node[anion] at (5.0,1.0) {\tiny Cl$^-$};
\node[anion] at (4.0,3.0) {\tiny Cl$^-$};
\node[anion] at (3.0,0.6) {\tiny Cl$^-$};
\node[anion] at (2.5,2.8) {\tiny Cl$^-$};
\node[anion] at (3.3,2.8) {\tiny Cl$^-$};
\node[anion] at (6.0,1.0) {\tiny Cl$^-$};
\draw[->,red] (5.5,2.5) -- (5.0,2.5) node[midway,above] {\tiny Cl$^-$};
\draw[->,red] (2.15,2.8) -- (1.7,2.8) node[midway,above] {\tiny Cl$^-$};

\node[zincion] at (1.39,1.5) {\tiny Zn$^{2+}$};
\node[zincion] at (1.6,2.2) {\tiny Zn$^{2+}$};
\node[zincion] at (1.4,0.8) {\tiny Zn$^{2+}$};

\draw[fill=white] (7.3,3.0) circle (0.2);
\draw[fill=white] (7.1,2.6) circle (0.2);
\node at (7.3,3.0) {\tiny H$_2$};
\node at (7.1,2.6) {\tiny H$_2$};

\node at (1.3,0.25) {\tiny Zn $\rightarrow$ Zn$^{2+}$ + 2e$^{-}$};
\node at (6.1,0.25) {\tiny 2H$_2$O + 2e$^{-}$ $\rightarrow$ H$_2$ + 2OH$^{-}$};

\draw[thick] (0.75,3.5) -- (0.75,5) -- (7.25,5) -- (7.25,3.5);
\draw[->,very thick] (1.3,4.9) -- (2.3,4.9) node[midway,above] {\scriptsize e$^{-}$ flow};
\node[draw,fill=yellow!20] at (4,5) {\small Voltmeter};

\end{tikzpicture}
\caption{{Schematic of a simple electrochemical cell with zinc and copper electrodes in NaCl electrolyte solution.}}
\label{Fig:1} 
\end{figure}

\vspace{0.2cm}

As far as computational studies are concerned, there is hardly any literature on the numerical approximations for modified PNP/NS  system for compressible electrolytes. While extensive research has been conducted on the mathematical and numerical investigation of classical PNP equations \cite{Prohl2009, Jerome1991, Yang2013} and the classical PNP/NS coupling system \cite{Ray2012, Prohl2010, Schmuck2011, Schmuck2009}. Nonetheless, the standard Finite Element Method (FEM) has been the most frequently employed approach in recent years. A different variant of the modified PNP/NS  system is explored in \cite{Mingyan}, where the standard Poisson equation in the PNP system is replaced with a fourth-order elliptic equation to achieve a more accurate representation of the electrostatic potential. A fully nonlinear Crank–Nicolson FEM for PNP equations was developed by Sun et al. \cite{Sun2016}, while Gao and He \cite{Gao2017} established a linearized finite element discretization and obtained unconditionally optimal error estimates for all governing variables. 
Beyond these developments, significant progress has been made in structure-preserving and energy-stable discretizations for electrokinetic systems, particularly using finite element and discontinuous Galerkin (DG) methods to maintain essential physical constraints such as positivity, charge conservation, and energy dissipation \cite{Cao2025, fu2022high, yang2024unconditional, liu2022positivity}. However, due to the coupling of different mechanisms, establishing the well-posedness of the coupled PNP/NS equations remains a challenging task. The existence and uniqueness of solutions for steady-state PNP equations have been addressed in \cite{Jerome1985, Park1997}, while similar results for the classical PNP/NS coupling system are discussed in \cite{Schmuck2009, Jerome2002}. Additionally, several other methods, including discontinuous Galerkin, mixed, conservative, stabilized, and weak Galerkin approaches, have been developed for both classical PNP and classical PNP/NS equations \cite{Gao2018, He2017, He2018, Kim2022, Linga2020, Xie2020, Prohl2009}.

\vspace{0.2cm}
\subsection{\large \textbf{Contributions}}
Reliable computational results may be challenging to obtain for modified PNP/NS coupling system as the diffusion fluxes $\bm{J}_{\alpha}$ are implicitly encoded. Additionally, unlike the non-dimensional form of the classical PNP/NS coupling system \cite{Druz2013, Groos2019}, the non-dimensional form of the modified PNP/NS coupling system is characterized by fractional non-linearity (refer to section \ref{Sec:3}), leading to singularities in the governing variables. Moreover, the presence of double layers in the electrical fields in the proximity of the liquid–solid interface mandates a fine spatio-temporal resolution.

\vspace{0.2cm}
Given these notable challenges, our objective is to numerically evaluate the validity of replacing the NP equation in the PNP/NS model with the formulations presented in  \eqref{Eq19nn} and \eqref{Eq8nn}.
 To achieve this aim, we employ the finite element method for simulating the modified PNP/NS model (under the mechanical equilibrium) across diverse spatial dimensions while considering both compressible and incompressible variants of fluid. A key contribution of this work is the reduction of the equilibrium system to a modified Poisson–Boltzmann formulation, providing additional analytical insight into the compressible case. To emphasize practical relevance, we also discuss the implications of compressible electrolytes in the context of double-layer capacitance behavior.
We further replicate a variety of operating conditions, including temperature fluctuations and variations in the bulk modulus, which influence the fluid’s compressibility and incompressibility.
Numerical simulations over different domains, including rectangular and annular battery designs, confirm the accuracy and robustness of the proposed computational scheme. The results also demonstrate that the limiting behavior in the incompressible regime aligns with the theoretical trends anticipated by Dreyer et al.~\cite{dreyer2013overcoming}.

\subsection{\large \textbf{Outline}}
 The remainder of this paper is organized into five sections. In Section \ref{Sec:2}, we present the governing equations of our model. 
Section \ref{Sec:3} provides the mathematical reformulation of the model in thermodynamic equilibrium, the associated dimensionless system and boundary conditions. Section \ref{Sec:4} is designated for the weak formulation related to the proposed model. In Section \ref{Sec:5}, we present the practical application of our research through numerical simulations of the modified PNP/NS model, including the parametric studies. We conclude our work in Section \ref{Sec:6} by noting its main contributions and proposing further research avenues in this domain.


\section{The model problem}
\label{Sec:2}

Let $\Omega \subset \mathbb{R}^n$ $(n\leq3)$ be a bounded Lipschitz domain. For given final time $t_F>0$ and a mixture with $N-1$ ionic species and a neutral solvent, the full set of equations for the modified PNP/NS system \cite{dreyer2013overcoming} in $\Omega \times [0,t_F]$ is given by 
\begin{subequations}
\label{Eq23n}
\begin{align}
&\frac{\partial (m_{\alpha}n_{\alpha})}{\partial t}  + \nabla\cdot(m_{\alpha}n_{\alpha}\bm{v}+\bm{J}_{\alpha})=0, \quad \alpha\;\in\;\{1,\ldots,N-1\}, \label{Eq19n}\\
&\frac{\partial \rho}{\partial t}  + \nabla\cdot(\rho\bm{v})=0,\label{Eq20n}\\
&\frac{\partial (\rho\bm{v})}{\partial t} + \nabla\cdot(\rho\bm{v}	\otimes\bm{v})+ \nabla p=-(\sum_{\alpha=1}^{N}z_{\alpha}e_0n_{\alpha})\nabla\varphi,\label{Eq21n}\\
&-\epsilon_0\epsilon_r\Delta\varphi =  \sum_{\alpha=1}^{N}z_{\alpha}e_0n_{\alpha}, \label{Eq22n}
\end{align}
\end{subequations} 
where, $\partial_t=\partial/\partial t$, $n_{\alpha}$ ($\alpha=1, 2, \ldots N$) are the number densities of ionic species with mass $m_{\alpha}$ and charge $z_{\alpha}$; $\varphi$ is electrostatic potential,  $\bm{v}$, $\rho$ and $p$ represent the velocity, total mass density and pressure of the fluid, respectively; and $\epsilon_0, \epsilon_r$ and $e_0$ respectively denote the dielectric constant, relative dielectric permittivity and elementary charge (fixed positive constants).

For a given temperature $T$ and positive definite kinetic matrix  $H_{\alpha\beta}$, the diffusion fluxes $\bm{J}_{\alpha}$  in equation \eqref{Eq19n} can be expressed as follows  
\begin{equation}
	\label{Eq8n}
	 \bm{J}_{\alpha}= - \sum_{\beta=1}^{N-1}H_{\alpha\beta} \bigg(\frac{\nabla\mu_{\beta}-\nabla\mu_N}{T}+ \frac{e_0}{T}\bigg(\frac{z_{\beta}}{m_{\beta}}-\frac{z_{N}}{m_{N}}\bigg)\nabla\varphi\bigg), \quad\alpha \in \{1, 2, \ldots N-1\}, 
\end{equation}
where
\begin{equation}
\label{Eq9n}
\mu_{\alpha}= g_{\alpha}^R+\frac{kT}{m_{\alpha}}\mbox{ln}\bigg(\frac{n_\alpha}{n}\bigg)+\frac{K}{m_{\alpha}n^R}\mbox{ln}\bigg(1+\frac{p-p^R}{K}\bigg) , \quad\quad\alpha \in \{1, 2, \ldots N\}. 
\end{equation}
\begin{equation}
\label{Eq10n}
p= p^R+K\bigg(\frac{n}{n^R}-1\bigg), 
\end{equation}
such that chemical potentials $\mu_{\alpha}$ $(\alpha=1,2\ldots N)$ are intricately related to number densities and specific Gibbs energies $(g_{\alpha})$. Here, $R$ indicates the reference state; $k$ and $K$ denote the Boltzmann constant and bulk modulus, respectively, and $n$ denotes the total number density, which is the sum of all the $n_{\alpha}$ $(\alpha=1,2\ldots N)$. Further details on the mathematical formulation of the modified PNP/NS coupling system are provided in \ref{Sec:2app}, and the notations used in this article are defined in \ref{Notations}. The boundary conditions related to the model will be discussed in the next section.


\section{The thermodynamical equilibrium based model problem in non-dimensional form}
\label{Sec:3}
We consider the mixture in thermodynamical equilibrium. 
As described in \cite{dreyer2013overcoming}, a system is said to be in thermodynamical equilibrium when velocity field $\bm{v}=0$ and diffusive fluxes $\bm{J}_{\alpha}=0$ for $\alpha = 1, 2, \ldots N-1$. Since the solvent is electrically neutral, we set $z_N = 0$. Further, assuming that the system is time-independent, the governing equations \eqref{Eq23n} reduce to:
\begin{subequations}
\label{Eq23*n}
\begin{align}
&\nabla\cdot(\bm{J}_{\alpha})=0, \quad \alpha\;\in\;\{1,\ldots,N-1\}, \label{Eq19*n}\\
& \nabla p=-(\sum_{\alpha=1}^{N-1}z_{\alpha}e_0n_{\alpha})\nabla\varphi,\label{Eq21*n}\\
&-\epsilon_0\epsilon_r\Delta\varphi =  \sum_{\alpha=1}^{N-1}z_{\alpha}e_0n_{\alpha}. \label{Eq22*n}
\end{align}
\end{subequations} 
As there are less number of equations than unknowns in the above system, we need to use constitutive laws \eqref{Eq9n} and \eqref{Eq10n} in diffusive fluxes \eqref{Eq8n}, which will subsequently be substituted into equation \eqref{Eq19*n}. Moreover, equation \eqref{Eq9n} suggests transforming the variable from  number densities to atomic fractions $(y_{\alpha})_{\alpha\in \{1, 2\ldots N\}}$ such that 
 \begin{equation}
	\label{Eq24}
y_\alpha=\frac{n_\alpha}{n}. 
\end{equation}
Given that $n= \sum_{\alpha=1}^{N}n_{\alpha}$, we have the following constraint 
 \begin{equation}
 \label{Eq25}
\sum_{\alpha=1}^{N}y_\alpha=1.
\end{equation}
Now, we aim to reduce the model \eqref{Eq19*n} to a closed system for the variables $(y_1, y_2,\ldots y_{N-1}, n, \varphi)$ by systematically utilizing the expressions \eqref{Eq8n}-\eqref{Eq10n}. To achieve this, we proceed with the following steps:

\vspace{0.5em}
\noindent \textbf{{Step 1: Simplifying the momentum and potential equations}}:

Using equations: \eqref{Eq10n} and \eqref{Eq24} in Equations \eqref{Eq21*n} and \eqref{Eq22*n}, we get
\begin{subequations}
\begin{align}
\frac{K}{n^R} \nabla n&=-\sum_{\alpha=1}^{N-1}e_0(z_\alpha y_\alpha)n\nabla\varphi,\label{Eq33}\\
-\epsilon_0\epsilon_r\Delta\varphi &= \sum_{\alpha=1}^{N-1}e_0(z_\alpha y_\alpha)n\label{Eq34}.
\end{align}
\end{subequations}
\vspace{0.5em}
\noindent \textbf{{Step 2: Simplifying the ionic concentration equations or reformulating the diffusive fluxes}}:

By substituting equation \eqref{Eq10n} into \eqref{Eq9n}, and employing the transformation \eqref{Eq24}, it is easy to see that 
\begin{align}
    \label{Eq9nnn}
\nabla\mu_{\alpha}= \frac{kT}{m_{\alpha}}\left(\frac{1}{y_{\alpha}}\nabla y_{\alpha}\right)+\frac{K}{n^Rm_{\alpha}}\left(\nabla \Tilde{n}\right) , \quad\quad\alpha \in \{1, 2, \ldots N\},
\end{align}
where
 \begin{equation}
 \label{Eq29}
\nabla \Tilde{n}= \frac{1}{n}\nabla n.
\end{equation}
Substituting the above equation into diffusive flux \eqref{Eq8n}, we obtain 
\begin{align}
	\label{Eq26}
	 \bm{J}_{\alpha}= &- \sum_{\beta=1}^{N-1}\frac{H_{\alpha\beta}K}{n^RT} \bigg(\frac{1}{m_{\beta}}-\frac{1}{m_{N}}\bigg)\nabla \Tilde{n}-
  \sum_{\beta=1}^{N-1}{H_{\alpha\beta}k}\bigg(\frac{1}{m_{\beta}y_{\beta}}\nabla y_{\beta}-\frac{1}{m_{N}y_{N}}\nabla y_{N}\bigg)\nonumber\\
  &-\sum_{\beta=1}^{N-1}\frac{H_{\alpha\beta}e_0}{T}\bigg(\frac{z_{\beta}}{m_{\beta}}-\frac{z_{N}}{m_{N}}\bigg)\nabla\varphi, \quad\quad\quad\quad\alpha \in \{1, 2, \ldots N-1\}.
\end{align}
To further simplify the system, we assume that the kinetic matrix $H_{\alpha\beta}$ is diagonal. This
assumption ideally approximates many mixtures \cite{dreyer2013overcoming}. Thus, we have
 \begin{equation}
 \label{Eq27}
H_{\alpha\beta}=H_{\alpha}\delta_{\alpha\beta}.
\end{equation}
Using equations \eqref{Eq25} and \eqref{Eq27} in \eqref{Eq26} with a neutral solvent ($z_N=0$) to eliminate $y_N$, we obtain
\begin{align}
	\label{Eq28}
	 \bm{J}_{\alpha}= &- \frac{H_{\alpha}K}{n^RT} \bigg(\frac{1}{m_{\alpha}}-\frac{1}{m_{N}}\bigg)\nabla \Tilde{n}-
  {H_{\alpha}k}\bigg(\frac{1}{m_{\alpha}y_{\alpha}}\nabla y_{\alpha}+\frac{1}{m_{N}y_{N}}\big(\sum_{i=1}^{N-1}\nabla y_{i}\big)\bigg)\nonumber\\
  &-\frac{H_{\alpha}e_0}{T}\frac{z_{\alpha}}{m_{\alpha}}\nabla\varphi, \quad\quad\quad\alpha \in \{1, 2, \ldots N-1\}. 
\end{align}
\noindent\textbf{{Step 3: Reformulating the diffusive flux in terms of target variables:}}  

The flux expression in \eqref{Eq28} is already written in terms of the key variables \( (y_1, \ldots, y_{N-1}, n, \varphi) \). However, to bring it closer to the structure of the classical Nernst–Planck model (cf. equation \eqref{Eq8nnn}), we utilize equations \eqref{Eq33} and \eqref{Eq29}, which yield:
 \begin{equation}
 \label{Eq31}
\nabla \Tilde{n} = -\frac{n^Re_0}{K}\sum_{i=1}^{N-1}(z_{i}y_{i})\nabla \varphi.
\end{equation}
Substituting equation \eqref{Eq31} into \eqref{Eq28}, we arrive at a flux expression that is fully represented in terms of \( (y_1, \ldots, y_{N-1}, \varphi) \):
\begin{align}
\label{Eq32}
\bm{J}_{\alpha} = 
& - H_{\alpha}k \left( \frac{1}{m_{\alpha} y_{\alpha}} \nabla y_{\alpha} + \frac{1}{m_{N} y_{N}} \sum_{i=1}^{N-1} \nabla y_{i} \right) \nonumber \\
& + \left[ \frac{H_{\alpha} e_0}{T} \left( \frac{1}{m_{\alpha}} - \frac{1}{m_{N}} \right) \sum_{i=1}^{N-1} z_i y_i 
- \frac{H_{\alpha} e_0}{T} \frac{z_{\alpha}}{m_{\alpha}} \right] \nabla \varphi,
\quad \alpha \in \{1, \ldots, N-1\}.
\end{align}
Introducing the coefficients \( D_{\alpha\alpha} \), \( D_{\alpha i} \), and \( D_{\alpha \varphi} \), the flux \eqref{Eq32} can be written compactly as:
{\begin{equation}
\label{Eq36}
\bm{J}_{\alpha} = D_{\alpha\alpha} \nabla y_{\alpha} 
+ \sum_{i=1}^{N-1} D_{\alpha i} \nabla y_i 
+ D_{\alpha \varphi} \nabla \varphi,
\quad \text{with } i \neq \alpha,\ \alpha \in \{1, \ldots, N-1\},
\end{equation}}
where the nonlinear coefficients are defined by:
\begin{subequations}
\begin{align}
D_{\alpha\alpha} &= - H_{\alpha}k \left( \frac{1}{m_{\alpha} y_{\alpha}} + \frac{1}{m_N y_N} \right), \quad \alpha \in \{1, \ldots, N-1\}, \label{Eq37} \\
D_{\alpha i} &= - H_{\alpha}k \cdot \frac{1}{m_N y_N}, \quad i \neq \alpha, \label{Eq38} \\
D_{\alpha\varphi} &= \frac{H_{\alpha} e_0}{T} \left[ \left( \frac{1}{m_{\alpha}} - \frac{1}{m_N} \right) \sum_{i=1}^{N-1} z_i y_i 
- \frac{z_{\alpha}}{m_{\alpha}} \right] , \quad \alpha \in \{1, \ldots, N-1\}. \label{Eq39}
\end{align}
\end{subequations} 
Thus, the model reduces to finding the variables $(y_1, y_2,\ldots y_{N-1}, n, \varphi)$. Once these variables are determined, the pressure for a compressible mixture can be readily calculated using \eqref{Eq10n}. However, in the case of an incompressible mixture, where the bulk modulus 
$K$ approaches infinity, the relation given by \eqref{Eq10n} is no longer valid. This particular case will be discussed in the latter portion of the article.

\subsection{Dimensionless form of modified PNP/NS system} \label{dimensionless}
The non-dimensional form of the classical PNP/NS coupling system is thoroughly discussed in \cite{Druz2013, Groos2019}. In recognition of the importance of numerical simulations, we similarly introduce dimensionless variables for the modified PNP/NS coupling system. These dimensionless quantities not only scale physical values to reference levels but also facilitate a more robust and stable computational environment. By removing the dependence on specific units, these dimensionless parameters allow for a smooth transition between different scales, promoting numerical stability and efficiency across diverse scenarios. Therefore, we introduce the
dimensionless variables and express the mass fractions ${M}_\alpha$ for each of the constituents as
 \begin{equation}
 \label{Eq40}
{M}_\alpha= \frac{m_\alpha}{m_N}.
\end{equation}
Additionally, we introduce the following dimensionless transformations 
 \begin{equation}
 \varphi \rightarrow \frac{\varphi}{\varphi_{BC}},\;n \rightarrow \frac{n}{n_0},\; \nabla \rightarrow x_0\nabla.   \nonumber
\end{equation}
Here, $n_0, \varphi_{BC}$ and $x_0$ denote the average value of the total number density, characteristic boundary potential, and characteristic length scale of the system, respectively. Furthermore, we introduce the new dimensionless parameters $\Psi=\frac{e_0}{kT}\varphi_{BC}$, $\Lambda= \frac{e_0^2n_0x_0^2}{\epsilon_0kT}$ and $\hat{K}=\frac{K}{n_0kT}$. 

Utilizing these parameters, the dimensionless form of equation \eqref{Eq19*n} is $\nabla\cdot(\bm{J}_{\alpha})=0$, where $\bm{J}_{\alpha}$ is given by equation \eqref{Eq36}. The corresponding dimensionless coefficients appearing in \( \bm{J}_{\alpha} \) are expressed as:
\begin{subequations}
 \begin{align}
D_{\alpha\alpha} &= -\bigg(\frac{1}{M_{\alpha}y_{\alpha}}+\frac{1}{y_{N}}\bigg), \label{Eq41}\\
D_{\alpha i} &= -\frac{1}{y_{N}}, \quad i=1,2,\ldots N-1 \;\mbox{and}\; i\neq \alpha,\label{Eq42}\\
D_{\alpha\varphi} &=  \Psi \bigg(\frac{1}{M_{\alpha}}-1\bigg)\sum_{i=1}^{N-1}(z_{i}y_{i}) -\Psi \frac{z_{\alpha}}{M_{\alpha}}.\label{Eq43}
\end{align}
\end{subequations}
Furthermore, by choosing \( n^R = n^0 \), the momentum and potential equations~\eqref{Eq33} and~\eqref{Eq34}
 transform into
\begin{subequations}
\begin{align}
 \nabla n&=-\frac{\Psi}{\hat{K}}\sum_{\alpha=1}^{N-1}(z_\alpha y_\alpha)n\nabla\varphi,\label{Eq44}\\
-\Delta\varphi &= \frac{\Lambda}{\Psi}\frac{1}{\epsilon_r}\sum_{\alpha=1}^{N-1}(z_\alpha y_\alpha)n\label{Eq45}.
\end{align}
\end{subequations}
Finally, the dimensionless form of the thermodynamical equilibrium based modified PNP/NS system is given as 
{\begin{subequations}
	\begin{align}
	\nabla \cdot& \bigg(
    D_{\alpha\alpha} \nabla y_{\alpha} 
    + D_{\alpha\varphi} \nabla \varphi 
    + \sum_{i=1}^{N-1} D_{\alpha i} \nabla y_i
\bigg) = 0,
\quad \text{with } i \neq \alpha,\quad \alpha \in \{1, \ldots, N-1\},\label{Eq44aa}\\
	\nabla n&=-\frac{\Psi}{\hat{K}}\sum_{\alpha=1}^{N-1}(z_\alpha y_\alpha)n\nabla\varphi,\label{Eq44n}\\
	-\Delta\varphi &= \frac{\Lambda}{\Psi}\frac{1}{\epsilon_r}\sum_{\alpha=1}^{N-1}(z_\alpha y_\alpha)n\label{Eq45n},
	\end{align}
\end{subequations}}
where, the nonlinear dimensionless coefficients $D_{\alpha\alpha}, D_{\alpha i}$ and $D_{\alpha\varphi}$ are given by equations \eqref{Eq41}-\eqref{Eq43}, respectively.

\subsection{Boundary conditions} \label{Bdry}

Assume that $\partial\Omega$ denotes the boundary of domain $\Omega$ such that $\partial\Omega=\Gamma_D\cup\Gamma_N$ with $ \mathring{\Gamma}_D \cap  \mathring{\Gamma}_N = \emptyset$ and  ${\left|\Omega\right|}$ be the measure of $\Omega$.   For the electric potential $\varphi$, we prescribe the following boundary conditions 
\begin{align}
\label{eq49}
\left\{
\begin{array}{ll}
\varphi=g & \mbox{on} \;\Gamma_D, \\
\frac{\partial \varphi}{\partial \vec{n}}=0 &\;\mbox{on}\;\Gamma_N.
\end{array}
\right.  
\end{align}
We assume no flux boundary conditions to ensure that there is no net flux of any ionic species through the boundary. This can be mathematically expressed as
\begin{equation}
{\bm{J}_{\alpha}}\cdot\vec{n}=0, \quad \quad {\alpha \in \{1, 2, \ldots N-1\}}. \label{eq51}
\end{equation}
Taking the divergence of \eqref{Eq44n}, we obtain 
\begin{equation}
\nabla\cdot\bigg(\nabla n+\frac{\Psi}{\hat{K}}\sum_{\alpha=1}^{N-1}(z_\alpha y_\alpha)n\nabla\varphi\bigg)=0.\label{eq51n}
\end{equation}
Again, multiplying \eqref{Eq44n} with the normal vector $\vec{n}$, we have
\begin{equation}
\bigg(\nabla n+\frac{\Psi}{\hat{K}}\sum_{\alpha=1}^{N-1}(z_\alpha y_\alpha)n\nabla\varphi\bigg)\cdot \vec{n}=0,\label{eq51nn}
\end{equation}
which serves as the Neumann boundary condition for total number density $n$. Equation \eqref{eq51n} alone does not guarantee that \eqref{Eq44n} holds due to the presence of an integration constant. However, by setting \eqref{eq51nn} as the Neumann boundary condition, the constant is eliminated, and thus, equations \eqref{Eq44n} and \eqref{eq51n} become equivalent when combined with the Neumann boundary condition \eqref{eq51nn}. This step brings \eqref{Eq44n} into alignment with the other equations, producing a second-order scaler equation and enabling a symmetric variational formulation that can be efficiently handled by standard FEM tools such as FEniCS \cite{logg2012finite}.

Further, the following side conditions arise 
\begin{align}
\stackinset{c}{}{c}{}{\displaystyle\int}{\rule{1.1ex}{0.2ex}}_{\Omega} y_{\alpha} \, dx&= y_{\alpha}^0 ,\quad \alpha=1,2 \ldots N-1,\label{eq50}\\ 
\stackinset{c}{}{c}{}{\displaystyle\int}{\rule{1.1ex}{0.2ex}}_{\Omega} n \, dx&= n^0,\label{eq52}
\end{align}
where, \;$\stackinset{c}{}{c}{}{\displaystyle\int}{\rule{1.1ex}{0.2ex}}_{\Omega} U(x) \, dx = \frac{1}{\left|\Omega\right|}\bigintss_{\Omega}U(x)dx.$
In case of Neumann boundary conditions, the above side conditions serve as integral constraints that eliminate the indeterminacy associated with additive constants, thereby ensuring the numerical uniqueness of the solution \cite{kergrene2017approximation, Andreianov}. Although the continuous problem may admit a family of solutions differing by a constant (unique upto a constant), our numerical scheme guarantees a unique discrete solution by penalizing components in the null space. These constraints are implemented numerically through Lagrange multipliers in the discretized system, as detailed in the next section. 
\subsection{{Modified Poisson-Boltzmann system}} Before presenting the weak formulation of the model, we first derive the modified Poisson–Boltzmann reduction in a general setting. The corresponding reduction for the incompressible case has already been carried out in \cite{dreyer2013overcoming}. Here, we introduce the isothermal assumption and recall that the system is in equilibrium. Referring to equation (70) in \cite{dreyer2013overcoming} and using equation \eqref{Eq28} above, we obtain
\begin{equation}\label{mm1}
\mathbf{J}_{\alpha}
= -\frac{H_{\alpha}}{T m_{\alpha}}\,
\nabla\!\Bigg[ kT\ln y_{\alpha}
+ \frac{K}{n^{R}}\ln\!\left(\frac{n}{n^{R}}\right)
+ e_0 z_{\alpha}\varphi \Bigg] = \bm{0}, 
\qquad \alpha=1,\dots,N-1.
\end{equation}

\medskip
\noindent\textbf{Dimensionless form:}  
Dividing equation \eqref{mm1} by \(kT\) leads to the dimensionless representation
\begin{equation} \label{mm2}
\frac{H_{\alpha}k}{ M_{\alpha}m_nx_0}\,
\nabla\Big( \ln y_{\alpha} + \hat{K}\ln n + z_{\alpha}\Psi\varphi \Big) = \bm{0}, \qquad \alpha=1,\dots,N-1.
\end{equation}
Therefore, the expression inside the parentheses must be constant in space:
\begin{equation}\label{mm3}
\ln y_{\alpha} + \hat{K}\ln n + z_{\alpha}\Psi\varphi = C_{\alpha},
\qquad \alpha=1,\dots,N-1.
\end{equation}

\medskip
\noindent\textbf{Equilibrium distribution for $N-1$ ionic species:}  
From \eqref{mm3} we obtain the compressible Boltzmann-type law
\begin{equation}\label{mm4}
y_{\alpha} = A_{\alpha}\, n^{-\hat{K}} \, e^{-z_{\alpha}\Psi\varphi}, 
\qquad A_{\alpha} := e^{C_{\alpha}}.
\end{equation}  
The side condition \eqref{eq50} determines the normalization constants \(A_\alpha\). Substituting \eqref{mm4} into \eqref{eq50} yields
\begin{equation}\label{mm5}
A_{\alpha}
= y_{\alpha}^0
\left(\frac{1}{|\Omega|}\int_\Omega n^{-\hat{K}} \, e^{-z_{\alpha}\Psi\varphi}\,dx \right)^{-1}.
\end{equation}
Thus, the equilibrium distribution consistent with the mean constraint is
\begin{equation}\label{mm6}
y_{\alpha} 
= y_{\alpha}^0\, \left({\dfrac{1}{|\Omega|}\int_\Omega n^{-\hat{K}} e^{-z_{\alpha}\Psi\varphi}\,dx}\right)^{-1} \,
\left({n^{-\hat{K}}\, e^{-z_{\alpha}\Psi\varphi}}\right),
\qquad \alpha=1,\dots,N-1.
\end{equation}

\medskip
\noindent\textbf{Equilibrium distribution of the number density (n):}  
For the neutral component \(N\), \(\nabla \mu_N = 0\) (see equ. (70), \cite{dreyer2013overcoming}) implies:
\begin{equation}\label{mm7}
y_N = A_N\, n^{-\hat{K}}, \quad A_N = e^{C_{N}},\; \text{for some constant} \; C_N.
\end{equation}
From the definition of atomic fractions,
\begin{equation}
\sum_{\alpha=1}^{N} y_\alpha = 1 
\quad \Rightarrow \quad 
\sum_{\alpha=1}^{N-1} y_\alpha + y_N = 1.
\end{equation}
Substituting the expressions for \(y_\alpha\) and \(y_N\) gives 
\begin{equation} \label{mm10}
\sum_{\alpha=1}^{N-1} A_\alpha\, n^{-\hat{K}} \, e^{-z_\alpha \Psi \varphi} + A_N\, n^{-\hat{K}} = 1 \quad \Rightarrow \quad n = \left( A_N + \sum_{\alpha=1}^{N-1} A_\alpha\, e^{-z_\alpha \Psi \varphi} \right)^{1/\hat{K}}.
\end{equation}
The constant \(A_N\) is determined from the total number density constraint \eqref{eq52}:
\begin{equation}\label{mm8}
\frac{1}{|\Omega|} \int_\Omega n \, dx = n^0
\quad \Rightarrow \quad
\frac{1}{|\Omega|} \int_\Omega \left( A_N + \sum_{\alpha=1}^{N-1} A_\alpha e^{-z_\alpha \Psi \varphi} \right)^{1/\hat{K}} dx = n^0.
\end{equation}
For general $\hat{K} \neq 1$, this is a nonlinear equation for $A_N$ and typically must be solved numerically. In the special case $\hat{K} = 1$, however, the equation becomes linear in $A_N$ and can be solved explicitly:
\begin{equation}\label{mm9}
A_N = n^0 - \frac{1}{|\Omega|} \sum_{\alpha=1}^{N-1} A_\alpha \int_\Omega e^{-z_\alpha \Psi \varphi} \, dx.
\end{equation}
Thus, the equilibrium distribution consistent with the mean density constraint for $\hat{K} = 1$ is 
\begin{equation}
n = n^0 - \frac{1}{|\Omega|} \sum_{\alpha=1}^{N-1} A_\alpha \int_\Omega e^{-z_\alpha \Psi \varphi} \, dx + \sum_{\alpha=1}^{N-1} A_\alpha \, e^{-z_\alpha \Psi \varphi}.
\end{equation}
Finally, substituting the above derived expressions into the governing equations \eqref{Eq45n} and \eqref{eq51n} yields a reduced formulation coupling $\varphi$ and $n$.

\section{Finite element approximation}
\label{Sec:4} This section discusses the finite element (FE) approximation including the associated weak formulation of the model. It covers the discretization of the computational domain, selection of finite elements and strategies for imposition of the boundary conditions. A fundamental understanding of these aspects provides the foundation for the subsequent discussion on the weak formulation.

\subsection{Gerneral description of the scheme}
\label{Sec:4.1} The origin of the Finite element method (FEM) lies in addressing the challenges posed by complex elasticity and structural analysis problems, particularly in the fields of civil and aeronautical engineering. One of the earliest contributions to the mathematical formulation of the FEM can be traced back to the works of Schellback \cite{Schellbach} and Courant \cite{courant1943variational}. One of the distinctive features of FEM is its ability to handle complex geometries through the use of mesh generation techniques. By dividing the computational domain into smaller subdomains or elements, FEM allows for adaptability in choosing the size and shape of these elements.

We begin by considering the finite element discretizations for the polyhedral domain $\Omega\subset\mathbb{R}^n, \;n\leq 3$. Let $\mathcal{T}_h$ denote a conformal partition of $\Omega$ into disjoint simplices $E$ of diameter $ h_E$ such that

  \begin{itemize}
  	\item  ${\overline \Omega = \cup_{E\in{\cal T}_h} E}$.
  	\item If $ E_1 $, $ E_2 $ are any two element of ${\cal T}_h $ such that $
  	E_1 \ne E_2 $, then either $E_1 \cap E_2  = \emptyset $ or $E_1 \cap E_2 $ is a common $r$-face of both simplices, for $0 \leq r \leq n-1$. Hence, for $n=2$, $E_1 \cap E_2 $ is empty or a common vertex or edge of both triangles.
  	
      \item  Let $h_E$ denote the length of longest edge ($1$-face length) of the simplex $E\in {\cal T}_{h}$ and define $h = \max_{E\in{\cal T}_{h}} h_{E}$ as the maximum edge length in the triangulation ${\cal T}_{h}$. We assume that for any fixed $h' > 0$, there exist two positive constants $K_1$ and $K_2$, such that for any triangulation $\mathcal{T}_h$ with $0 < h < h'$, the following condition holds:
  	\begin{eqnarray*}
  		K_1 h  \;\leq \;\mbox{diam}(E)\; \leq \;K_2 h, \quad \quad \forall\; E \in{\cal T}_h,
  	\end{eqnarray*}
        where $K_1$ and $K_2$ are independent of $h$.
  \end{itemize}

Based on the above triangulation ${\cal T}_{h}$, we consider the following finite-element space of continuous, piecewise linear functions:
 \begin{align*}
    S_h^1(\Omega, {\cal T}_h) &\coloneqq \left\{ u \in C^0(\overline{\Omega}) \;\; \big|\;\; u|_{E} \in \mathbb{P}_1, \;\; \forall \; E\in{\cal T}_h  \right\},
  \end{align*}
  having denoted by $\mathbb{P}_1$ the space of polynomials with degrees less than or equal to one. Note that the space $S_h^1(\Omega, {\cal T}_h)$ is a subspace of Sobolev space $H^1(\Omega)$. For a more comprehensive understanding of more general finite elements, interested readers may refer to the works of \cite{ankurKDV2023, jiwari2023new, ankur2023new, quarteroni2009numerical, prohl2010convergent}. We will use the above $\mathbb{P}_1$ finite element space to reduce our weak formulation on finite-dimensional space, which will be elaborated in the next part of this section.  

In equation \eqref{eq49}, if ${\Gamma}_D=\emptyset$ (pure Neumann on \(\partial\Omega\)), then the Poisson equation \eqref{Eq45n} becomes
 \begin{align}
    \begin{split}
        -\Delta\varphi &= n^F, \quad \text{in } \Omega, \quad \quad n^F := \frac{\Lambda}{\Psi}\frac{1}{\epsilon_r}\sum_{\alpha=1}^{N-1}(z_\alpha y_\alpha)n,  \\
        \frac{\partial \varphi}{\partial \Vec{n}} &= 0 \quad \text{on } \partial\Omega.
    \end{split}
    \label{eq57}
\end{align}
Upon integrating the above equations and employing the divergence theorem, it becomes evident that the following crucial compatibility condition must hold:
 \begin{align}
\int_{\Omega} n^F\, dx=  0. \label{eq56}
\end{align}
This quantity represents the total electrostatic charge stored in the system, which in general may not vanish (see Section \ref{Sec:6.7}). Consequently, the problem may be ill-posed.

On the other hand, as pure Neumann conditions are applied, $\varphi$
 is only determined up to a constant $c$. Thus, an additional constraint is imposed, for instance 
 \begin{align}
\int_{\Omega} \varphi\, dx&= 0.\label{eq55}
\end{align}
As outlined in \cite{kergrene2017approximation}, the problem \eqref{eq57} with the constraint \eqref{eq55} can be effectively addressed using the method of Lagrange multipliers. This approach yields a new well-posed variational problem such that the right-hand side $f$ is redefined with the help of a constant to satisfy the condition \eqref{eq56}. Similarly, the equations for $y_{\alpha}$ and $n$ ( \eqref{Eq44aa} and \eqref{eq51n} ) with integral constraints \eqref{eq50} and \eqref{eq52} are also enforced via Lagrange multipliers in the FEM implementation.


\subsection{Variational formulation}
\label{Sec:4.2}
In the rest of the article, we shall employ the standard mathematical notations for norms and Sobolev spaces. We utilize $\|\cdot\|_{m}$ for representing the norm in Sobolev spaces $H^m(\Omega)$ for $m\geq 1$. Also, in the case of Sobolev space $H^0(\Omega)$ or $L^2(\Omega)$, the inner product and the norm are denoted as  $(\cdot, \cdot)$ and  $\|\cdot \|$, respectively.
In addition, let us define the spaces
\begin{align}
&H(\Omega) = \left\{
  u \in H^1(\Omega) : \int_{\Omega} u \, dx = 0
\right\},\nonumber\\
&H^1(\Omega,\Gamma_D) = \left\{
  u \in H^1(\Omega) : u|_{\Gamma_D} = 0
\right\}.\nonumber
\end{align}
Before we proceed, one can notice that some of the governing variables utilize pure Neumann boundary conditions. Fredholm's alternative implies that those variables can be determined up to additive constants. For the uniqueness of solutions, side conditions \eqref{eq50} and \eqref{eq52} are very helpful. The conventional approach to imposing constraints involves the introduction of a Lagrange multiplier. Hence, proceeding in the usual way as done in \cite{babuvska2010residual} and following the notations of \cite{kergrene2017approximation}, a weak solution of the modified PNP/NS coupling system in thermodynamical equilibrium is a pair 
\begin{align}
    &(y_1, y_2, \ldots y_{N-1}, n, \varphi, c_1, c_2, \ldots , c_N)\nonumber\\&=(\Tilde{y}_1+y_1^0, \Tilde{y}_2+y_2^0, \ldots, \Tilde{y}_{N-1}+y_{N-1}^0,\Tilde{n}+n^0, \Tilde{\varphi}+g, c_1, c_2, \ldots ,c_N ) \in [H^1(\Omega)]^{N+1}\times \mathbb{R}^N\nonumber
\end{align}
such that $\Tilde{y}_{\alpha}\in H(\Omega)$ with $\alpha\in \{1,2,\ldots N-1\}$, $\Tilde{n}\in H(\Omega)$, $\Tilde{\varphi}\in H^1(\Omega,\Gamma_D)$ and 
\begin{align}
&\int_{\Omega} \bm{J}_{\alpha} \cdot \nabla v_{\alpha} \, dx - \int_{\partial\Omega} \bm{J}_{\alpha}\cdot \vec{n}v_{\alpha} \, ds +\int_{\Omega} c_{\alpha} v_{\alpha} \, dx=0 ,  \quad \alpha=1,2, \ldots, {N-1},\nonumber\\
&\int_{\Omega} \left( \nabla n + \frac{\Psi}{\hat{K}} \sum_{\alpha=1}^{N-1}(z_\alpha y_\alpha)n\nabla\varphi \right) \cdot \nabla v_n \, dx 
-\int_{\partial\Omega} \left( \nabla n + \frac{\Psi}{\hat{K}} \sum_{\alpha=1}^{N-1}(z_\alpha y_\alpha)n\nabla\varphi \right) \cdot \vec{n}v_n \, ds
\nonumber\\&\hspace{5cm}+\int_{\Omega} c_{n} v_{n} \, dx=0,\nonumber\\
& \int_{\Omega} \nabla\varphi \cdot \nabla v_\varphi \, dx - \int_{\Gamma_N}  \nabla\varphi \cdot \vec{n}v_\varphi \, ds - \int_{\Omega} \frac{\Lambda}{\Psi}\frac{1}{\epsilon_r}\sum_{\alpha=1}^{N-1}(z_\alpha y_\alpha)n v_\varphi \, dx
=0,\nonumber\\
& \int_{\Omega} y_{\alpha} d_{\alpha} \, dx -  \int_{\Omega} y_{\alpha}^0 d_{\alpha} \, dx=0, ,  \quad \alpha=1,2, \ldots, {N-1}, \nonumber\\
& \int_{\Omega} n d_{n} \, dx -  \int_{\Omega} n^0 d_{n} \, dx=0, 
\end{align}
for all $(v_1, v_2, \ldots v_{N-1}, v_n, v_\varphi, d_1, d_2, \ldots , d_N)\in [H^1(\Omega)]^{N}\times H^1(\Omega,\Gamma_D)\times \mathbb{R}^N$. Note that, given data $g\in H^{1/2}(\Omega)$ and $\Vec{n}$ denotes the outward unit normal. Now, we reduce the above weak formulation on the $\mathbb{P}_1$ finite element space, as defined in the last section, to effectively simulate our model. A comprehensive analysis, such as demonstrating well-posedness and deriving error estimates for modified PNP/NS coupled system, will be the subject of a follow-up paper.

\subsection{Implementation using FEniCS}
\label{Sec:4.3}
To implement the finite element framework, we leverage the capabilities of the open-source library FEniCS. It is a research and software project designed to facilitate the creation of software for automated computational mathematical modeling. Its core objective is to provide a platform for developing software that is easy to use, intuitive, efficient, and flexible for solving partial differential equations (PDEs) through FEM.
Originating in 2003, FEniCS has evolved through collaborative efforts among researchers from various universities and research institutes worldwide \cite{logg2012finite}. The library simplifies the computational process by automating the solution of PDEs, eliminating the need for manual intervention in certain steps of the FEM. In practice, using FEniCS typically involves supplying the weak formulation of a partial differential equation, accompanied by essential technical specifications such as function spaces, finite elements, and boundary conditions. FEniCS then undertakes the numerical solution process, automating steps like matrix assembly and equation system solving. This streamlined approach enhances the efficiency and accessibility of implementing Finite Element simulations, making it a valuable tool for researchers and practitioners in the field.

\section{Numerical results and discussion}
\label{Sec:5}
In this section, we simulate the present model using the Finite element scheme discussed in the previous section. We simulate the model for a wide range of parameter values, such as temperature and bulk modulus,  as they affect the numerical stability and qualitative behavior of the model. Throughout this section, we assume isothermal conditions, unless otherwise explicitly stated.
\subsection{\texorpdfstring{Empirical L\textsuperscript{2}-convergence}{Empirical L2-convergence}} \label{Sec:6.1} We first use proposed scheme for manufactured solutions to verify its robustness and accuracy. Let us assume that there are three constituents, namely cation (C), anion(A), and neutral solvent (S), having charges $z_C=1, z_A=-1, z_S=0$. Thus, the mixture can be described by field quantities such as electric potential $\varphi$, number density $n$, atomic fraction of cation $y_C$ and atomic fraction of anion $y_A$. Also, equation \eqref{Eq25} suggests choosing the manufactured solution such as $0< y_A, y_C< 1$. Let us propose the following exact solutions:
\begin{subequations}
\begin{align}
y_C&= \frac{1}{2+x^4+y^3},
&y_A = \frac{1}{3+x^3+y^2},\nonumber\\
\varphi &= \frac{1}{5+x^2+y^3}, 
&n = \frac{1}{1+x^3+y^3}, \nonumber
\end{align}
\end{subequations}
on the square domain $[0,1]^2$. Now, we impose the following non-homogeneous boundary conditions, specified in section \ref{Bdry}:
\begin{subequations}
\begin{align}
&~~~~~~~~~~~~~\vec{J_{C}}\cdot\vec{n}=g_1, \;\;  \vec{J_{A}}\cdot\vec{n}=g_2 \; \mbox{on}\; \partial\Omega,\nonumber \\
&~~~~~~~~~~~~~(\nabla n + \frac{\Psi}{\hat{K}}(z_Cy_C+z_Ay_A)n\nabla\varphi)\cdot\vec{n}=g_3\; \mbox{on}\; \partial\Omega,\nonumber
\end{align}\nonumber
\end{subequations}
$$\left\{
  \begin{array}{ll}
  \varphi=g_4, & \mbox{when}\; x=0 , \\
  \varphi=g_5, & \mbox{when}\; x=1 ,\\
  \frac{\partial \varphi}{\partial \Vec{n}}=g_6, &\mbox{elsewhere}\;\mbox{on}\;\partial\Omega.
  \end{array}
  \right.$$

\noindent Here, the functions $g_1, g_2$, $g_3, g_4$, $g_5$ and $g_6$ are determined according to the exact solutions, with the side conditions chosen to satisfy equations
\eqref{eq50}-\eqref{eq52}. 
 The constants are specified as $M _C=1$, $M _A=1$, $\hat{K}=1, \Lambda=1000$, $\Psi=1$ and $\chi=1$.
 \begin{table}[!h]
    \scriptsize
    \setlength{\tabcolsep}{3.5pt}
    \renewcommand{\arraystretch}{2}
    \caption{History of convergence and computational orders for linear finite elements in Example \ref{Sec:6.1}.}
    \centering
    \begin{tabular}{|c| c| c| c| c| c| c| c| c|}
        \hline
        & \multicolumn{4}{c|}{Cation} & \multicolumn{4}{c|}{Anion} \\  
        \hline
        \(h\) & \(L^2\) error & Order & \(H^1\) error & {Order} & \(L^2\) error & Order & {\(H^1\) error} & {Order} \\
        \hline
        \(1/4\) & \(1.7360 \times 10^{-3}\) & -- & \(3.6672 \times 10^{-2}\) & -- & \(1.2031 \times 10^{-3}\) & -- & \(1.4523 \times 10^{-2}\) & -- \\
        \(1/8\) & \(4.4202 \times 10^{-4}\) & \(1.9735\) & \(1.8708 \times 10^{-2}\) & \(0.9710\) & \(3.5898 \times 10^{-4}\) & \(1.7447\) & \(7.1546 \times 10^{-3}\) & \(1.0213\) \\
        \(1/16\) & \(1.1102 \times 10^{-4}\) & \(1.9932\) & \(9.4017 \times 10^{-3}\) & \(0.9926\) & \(9.3869 \times 10^{-5}\) & \(1.9352\) & \(3.5434 \times 10^{-3}\) & \(1.0137\) \\
        \(1/32\) & \(2.7786 \times 10^{-5}\) & \(1.9983\) & \(4.7077 \times 10^{-3}\) & \(0.9978\) & \(2.3738 \times 10^{-5}\) & \(1.9833\) & \(1.7668 \times 10^{-3}\) & \(1.0039\) \\
        \hline
    \end{tabular}
    \label{tbb1}
\end{table}

 \begin{table}[!h]
    \captionsetup{font=scriptsize}  
    \scriptsize
    \setlength{\tabcolsep}{3.5pt}
    \renewcommand{\arraystretch}{2}
    \caption{History of convergence and computational orders for linear finite elements in Example \ref{Sec:6.1}.}
    \centering
    \begin{tabular}{|c| c| c| c| c| c| c| c| c|}
        \hline
        & \multicolumn{4}{c|}{Potential} & \multicolumn{4}{c|}{Total Number Density} \\  
        \hline
        \(h\) & \(L^2\) error & Order & {\(H^1\) error} & {Order} 
        & \(L^2\) error & Order & {\(H^1\) error }& {Order} \\
        \hline
        \(1/4\) & \(2.8330 \times 10^{-3}\) & -- & \(1.5569 \times 10^{-2}\) & -- 
        & \(5.9769 \times 10^{-3}\) & -- & \(8.5963 \times 10^{-2}\) & -- \\
        
        \(1/8\) & \(8.2360 \times 10^{-4}\) & \(1.7823\) & \(5.4262 \times 10^{-3}\) & \(1.5207\) 
        & \(1.6260 \times 10^{-3}\) & \(1.8780\) & \(4.4692 \times 10^{-2}\) & \(0.9436\) \\
        
        \(1/16\) & \(2.1369 \times 10^{-4}\) & \(1.9464\) & \(1.9803 \times 10^{-3}\) & \(1.4542\) 
        & \(4.1547 \times 10^{-4}\) & \(1.9685\) & \(2.2597 \times 10^{-2}\) & \(0.9839\) \\
        
        \(1/32\) & \(5.3990 \times 10^{-5}\) & \(1.9847\) & \(8.5439 \times 10^{-4}\) & \(1.2127\) 
        & \(1.0443 \times 10^{-4}\) & \(1.9921\) & \(1.1333 \times 10^{-2}\) & \(0.9956\) \\
        \hline
    \end{tabular}
    \label{tbbb1}
\end{table}

\noindent We have employed the following relation to calculate the order of convergence 
$$\mbox{Order of convergence}= \frac{\log\big(\frac{\|e_{j+1}\|}{\|e_j\|}\big)}{\log\big(\frac{h_{j+1}}{h_j}\big)},$$
where $e_j$ denotes the error in $L^2$ norm on the jth iteration and $h_j$ is the corresponding mesh size.
\begin{figure}
	\centering
	\includegraphics[width=0.5\linewidth]{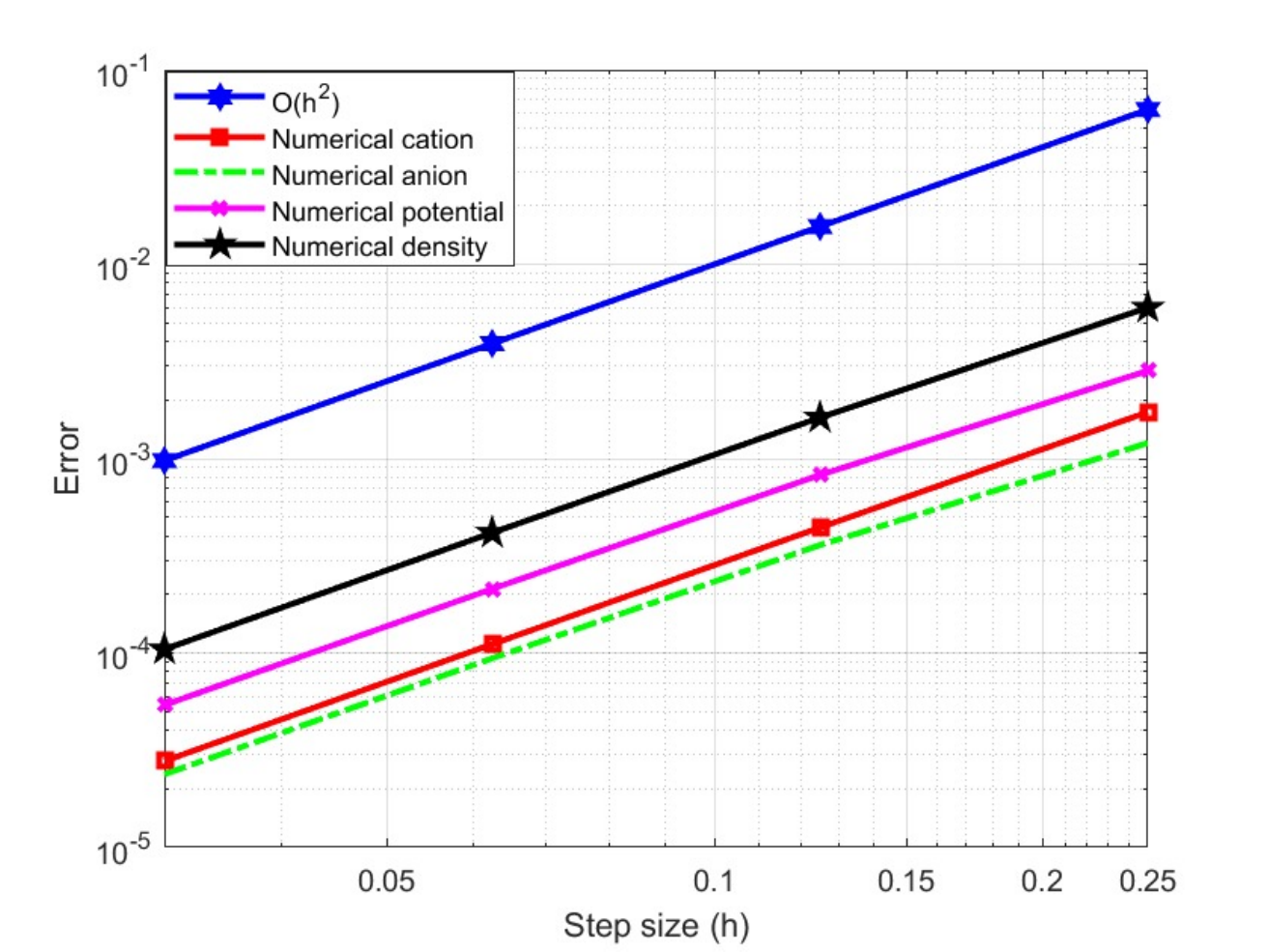}
	\caption{Log-Log plot showing second order of convergence 
          in \(L^{2}\) norm in example \ref{Sec:6.1}.}
	\label{Fig:2} 
\end{figure}

The convergence results for the space
discretization parameters $h = 0.25, 0.125, 0.0625$ and $0.03125$ are shown in Tables \ref{tbb1} and \ref{tbbb1}. They demonstrate second-order convergence in the \( L^2 \) norm and {first-order convergence in the \( H^1 \) norm for linear finite elements}, which aligns with the theoretical results observed in \cite{jiwari2023new, Jiwari2025-1, quarteroni2009numerical, ankur2024, Ismail2024}. A Log-Log plot, shown in Figure \eqref{Fig:2}, provides a visual representation of the error history in the $L^2$ norm. These comprehensive findings confirm the accuracy and reliability of our computational approach.

\begin{figure}
	\centering
	\includegraphics[width=0.7\linewidth]{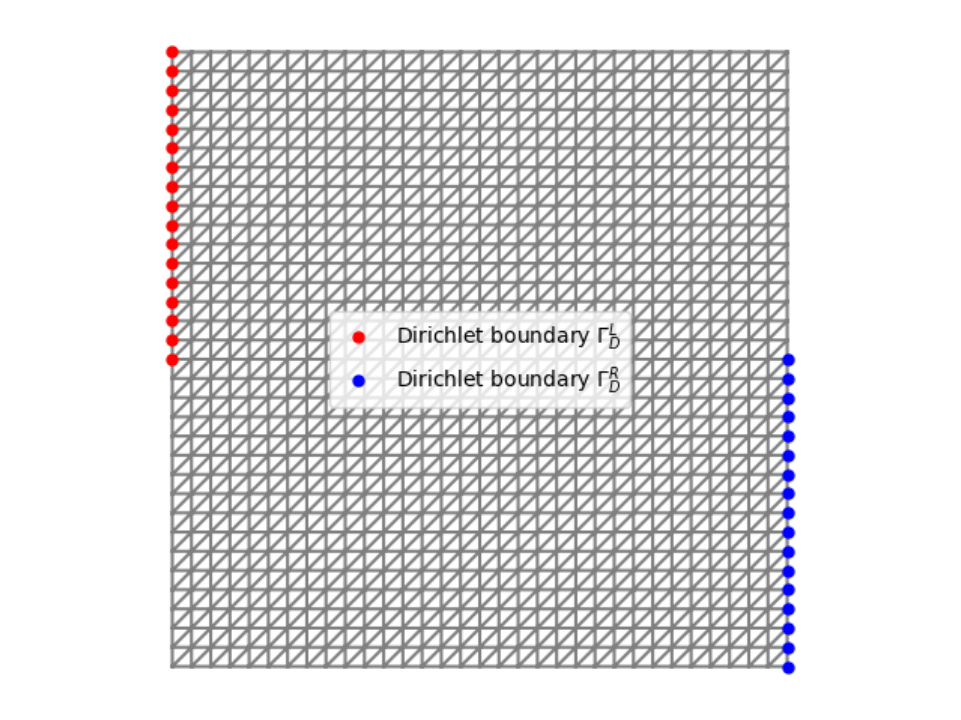}
	\caption{Boundary conditions on the electric potential \(\varphi\) for the unit square domain in Example~\ref{Sec:6.2}. All other boundary conditions for each variable are the same as those in Section~\ref{Bdry}.}
	\label{Fig:3} 
\end{figure}


\subsection{Symmetric 1:1 electrolyte\protect\footnote{An electrolyte having only one cationic species and one anionic species, both with same charge magnitude, for example, NaCl, CaS$O_4$.} dynamics in compressible regimes}
\label{Sec:6.2}

There are several ways to distinguish between incompressible and compressible mixtures. In this article, we define this concept with the help of Bulk modulus $K$ given by constitutive law \eqref{Eq10n}. Notably, the limit 
$K \rightarrow\infty$ signifies an incompressible mixture, a concept emphasized by Dreyer et al. \cite{dreyer2013overcoming}.
Since the dimensionless parameter $\Hat{K}$ is directly proportional to 
Bulk modulus $K$, we set $\Hat{K}=1$
Additional parameter settings are detailed in Table \ref{Compressible}. Using these parameters, we simulate the model on domain $[0,1]^n$ with $n\leq3$.  
We describe the electrodes through electrical potentials that are applied on the boundary of the domain. Assume that Dirichlet boundary $\Gamma_D=\Gamma_D^L\cup \Gamma_D^R$. We prescribe negative potential on the part of  Dirichlet boundary $\Gamma_D^L$ and positive potential on $\Gamma_D^R$ to set up a potential difference. For instance, Figure \eqref{Fig:3} illustrates the applied electric potential in the case of a unit square domain. Analogously, we apply electric potentials in 1D and 3D scenarios to govern ion transport phenomena within our mixture.

\begin{figure}[h!]
	\centering
	\includegraphics[width=0.7\linewidth]{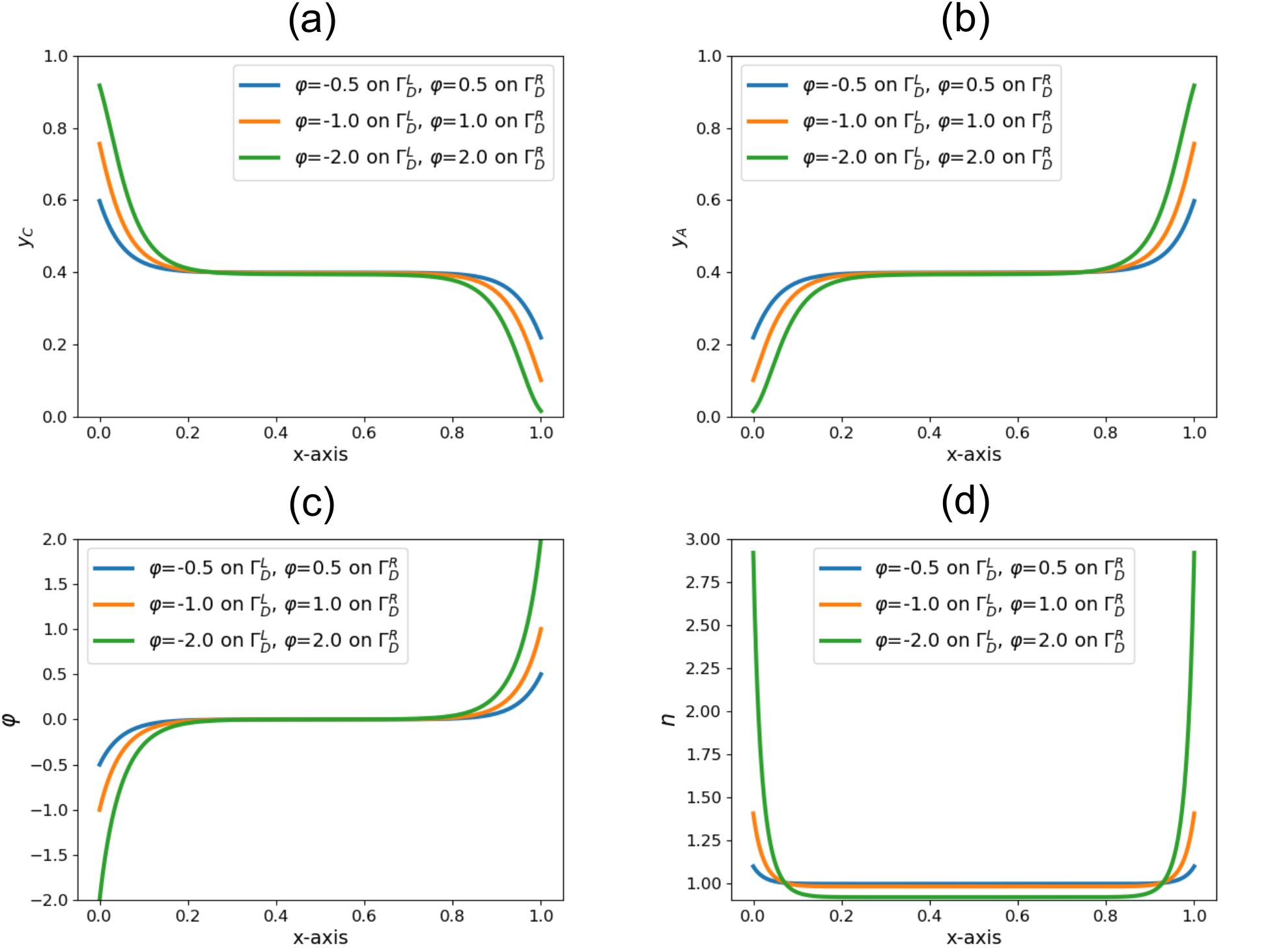}
	\caption{{Numerical solutions of (a) cation $y_c$, (b) anion $y_a$, (c) electric potential $\varphi$, and (d) total number density $n$ for different applied potentials on the unit interval in Example~\ref{Sec:6.2}. The parameter settings are detailed in Table~\ref{Compressible}, and the boundary conditions are described in Section~\ref{Bdry}}.}
	\label{Fig:4} 
\end{figure}

\begin{figure}[h!]
	\centering
	\includegraphics[width=0.7\linewidth]{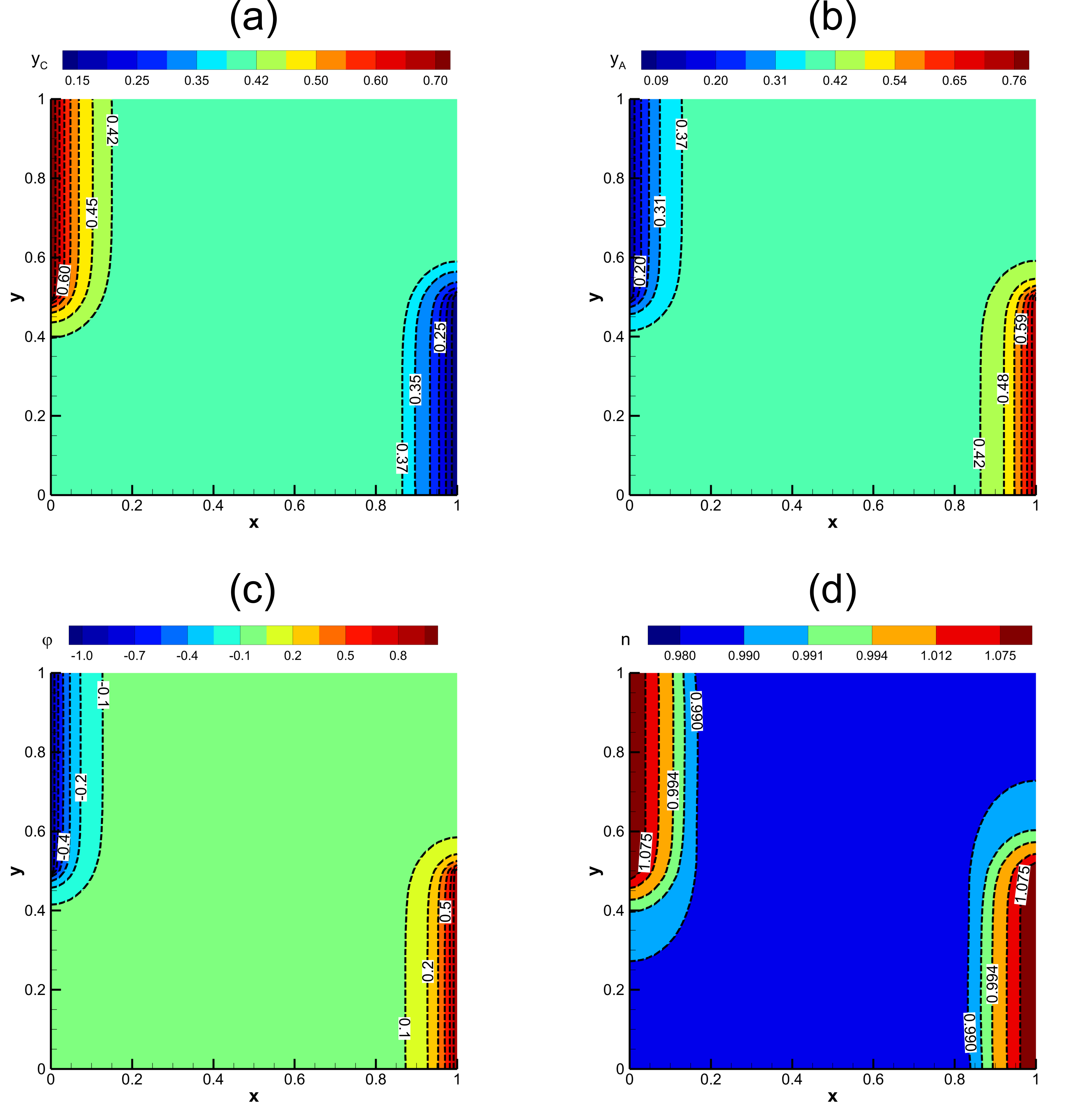}
	\caption{{Numerical solutions of (a) cation $y_c$, (b) anion $y_a$, (c) electric potential $\varphi$, and (d) total number density $n$ with applied potential $\pm 1$ unit on the 2D square domain in Example~\ref{Sec:6.2}. For parameter values, see Table~\ref{Compressible}; for boundary conditions, see Section~\ref{Bdry} and Fig.~\ref{Fig:3}.}}
	\label{Fig:5} 
\end{figure}

\begin{figure}[h!]
	\centering
	\includegraphics[width=0.7\linewidth]{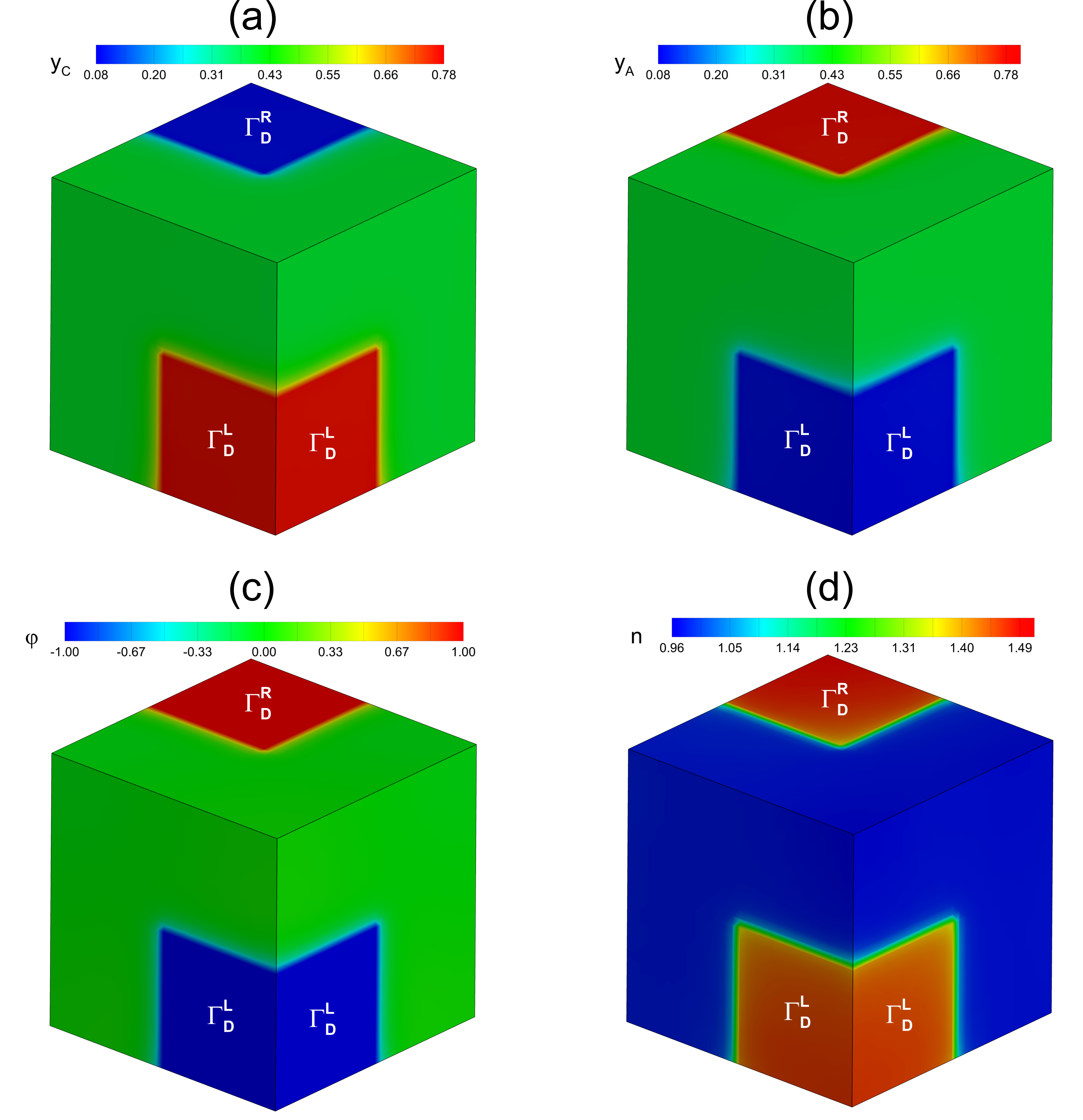}
	\caption{{Numerical solutions of (a) cation $y_c$, (b) anion $y_a$, (c) electric potential $\varphi$, and (d) total number density $n$ with applied potential $\pm 1$ unit for the 3D cubic domain in Example~\ref{Sec:6.2}. For parameters, see Table~\ref{Compressible}; for boundary conditions, see Section~\ref{Bdry}, and Fig.~\ref{Fig:3}, which is analogously extended from the unit square to the unit cube.}}
	\label{Fig:6} 
\end{figure}

\begin{table}[h!]
    \centering
    \caption{Parameters used for numerical calculation.}
    \label{Compressible}
    \footnotesize
    \begin{tabular}{lcc}
        \toprule
        \textbf{Description} & \textbf{Quantity} & \textbf{Value} \\
        \midrule
        Mass fraction of cations & $M_C$ & 0.1 \\
        Mass fraction of anions & $M_A$ & 0.1 \\
        Charge number of cations & $z_C$ & 1 \\
        Charge number of anions & $z_A$ & -1 \\
        Dielectric susceptibility & $\chi$ & 1 \\
        Average atomic fraction of cations & $y_{C}^0$ & 0.4 \\
        Average atomic fraction of anions & $y_{A}^0$ & 0.4\\
        Average total number density & $n^{0}$ & 1 \\
        Dimensionless parameter & $\Psi$ & 1 \\
        Dimensionless parameter & $\Lambda$ & 1000 \\
        Dimensionless parameter & $\hat{K}$ & 1 \\
        \bottomrule
    \end{tabular}
\end{table}
Numerical results for unit interval $[0,1]$ are shown in Figure \eqref{Fig:4}. We note that boundary layers are visible close to the electrodes in the concentration of all quantities: atomic fraction of cation ($y_C$), atomic fraction of anion ($y_A$), electric potential ($\varphi$), and number density ($n$). In the middle part of the domain, all of them have nearly a constant value equal to the chosen average value $y_C^0, y_A^0$ and $n^0$. We see that positively charged cations migrate towards the left electrode, subjected to a negative potential, while being repelled from the boundary where a positive potential is applied. Additionally, this concentration gradient along the electrodes intensifies with larger potential differences. We observe that negatively charged anions simply exhibit the opposite behavior. It is interesting to note that in the case of a compressible mixture, the number density is very high on both electrodes, showing that there is a significant accumulation of charge carriers near the electrodes. 
The numerical findings in the case of unit square domain \eqref{Fig:3} show that the concentration of the governing variable barely changes in the middle of the domain (See Figure \eqref{Fig:5}.  However, a distinct boundary layer appears along the Dirichlet borders $\Gamma_D^L$ and $\Gamma_D^R$, with a substantially greater number density. The atomic fraction of cations gets larger on the left side and lower on the right side and vice-versa for atomic fractions of anions. Additionally, the scalar field $\varphi$ also has a plateau value in the middle of the domain and changes in the boundary regions of the electrodes. It is smoothly distributed between the low and high Dirichlet values applied on both sides. 
The numerical results for the unit cube, where the potential difference is established in opposing corners of the cube, are presented in Figure \eqref{Fig:6}. It is noteworthy that these findings align with the observed trends in the unit square domain, indicating consistency in the numerical outcomes. Specifically, the concentration gradients persist prominently near the Dirichlet boundaries, while minimal alterations are noted in the central region of the domain.

\subsection{Symmetric 1:1 electrolyte: Transition to incompressible behavior with increasing Bulk modulus\label{Sec:6.3}} 
As earlier studies of this model have largely focused on the incompressible limit, our approach in this section explicitly accounts for spatial variations in the total ion concentration $n$ across different compressibility states.
According to Dreyer et al. \cite{dreyer2013overcoming}, the limit $K \rightarrow\infty$ signifies an incompressible mixture.  It is essential to note that, based on the constitutive law \eqref{Eq10n}, this limit does not imply that the pressure tends to infinity.  Instead, as $K \rightarrow\infty$, the number density $n$ approaches $ n^R= n^0$. However, the pressure can no longer be determined using constitutive law \eqref{Eq10n} as $K\left(\frac{n}{n^0}-1\right)$ becomes undetermined. 

\begin{figure}
	\centering
	\includegraphics[width=0.8\linewidth]{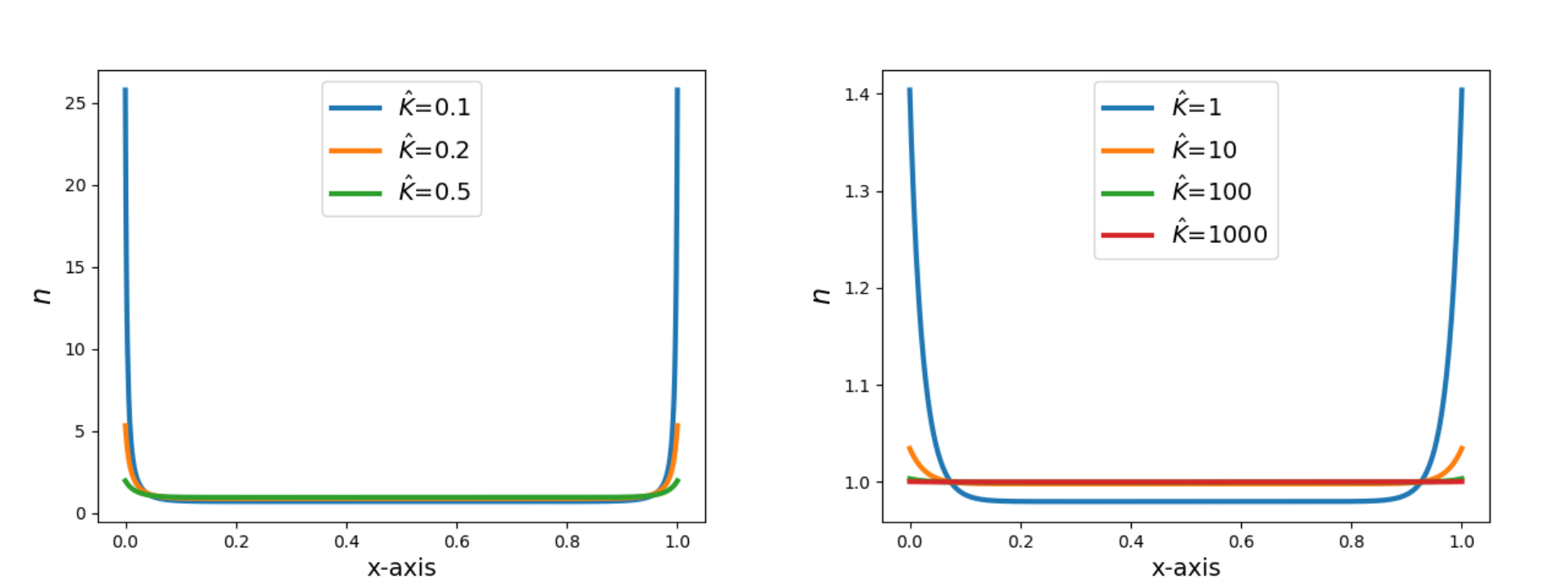}
	\caption{{Effect on number density over the unit interval with increasing bulk modulus in Example~\ref{Sec:6.3}. Parameter settings are same as Table~\ref{Compressible} except for varying values of $\hat{K}$. Boundary conditions are described in Section~\ref{Bdry}, with an applied potential of $-1$ unit at $x=0$ and $1$ unit at $x=1$.}}
	\label{Fig:7} 
\end{figure}

The most intuitive and significant validation of our compressible model is depicted in Figure \ref{Fig:7}. The simulated behavior of the number density aligns perfectly with the theoretical expectations described by Dreyer et al. \cite{dreyer2013overcoming}. This trend is observed across various values of $\Hat{K}$ ranging from $0.1$ to $1000$, with other parameters as per Table \ref{Compressible}. For a small bulk modulus ($\Hat{K}=0.1$), the electrolyte is highly compressible, leading to a significant accumulation of ions and thus a high local number density $n$ near the electrodes. As $\Hat{K}$ increases by orders of magnitude, the density profile flattens and $n$ converges uniformly towards the reference value $n^0$. This progression vividly illustrates the transition from a compressible to an incompressible state, where spatial variations in total concentration are suppressed. 

The physical implications of this density variation are profound and directly impact the electrolyte's ability to store charge. As discussed in Section \ref{Sec:6.7}, the enhanced ionic accommodation in compressible systems at high potentials, evident in Fig. \ref{Fig:7}, is the very mechanism that boosts the double-layer capacitance, $C_{dl}$. Conversely, as $K$ increases and the system becomes incompressible, the number density is constrained near $n^0$, reducing the space charge distribution and consequently lowering the overall $C_{dl}$. 
\subsection{{Ionic distributions in annular battery designs at equilibrium} } 
\label{Sec:6.4.1} 

{The annular geometry introduces unique characteristics to the design of battery systems. The electrode materials may be accommodated within the annular space, where the inner and outer surfaces function as the anode and cathode, respectively. This design provides a higher surface area for electrochemical reactions compared to a traditional rectangular cell. Moreover, the annular shape facilitates efficient heat dissipation, capitalizing on the expanded surface area available for thermal management. Such designs are particularly useful in applications where space is a critical factor, such as portable electronics or electric vehicles. Additionally, multiple annular cells can be combined in a modular fashion to create a battery pack with scalable capacity. }

\begin{figure}[h!]
	\centering
	\includegraphics[width=0.7\linewidth]{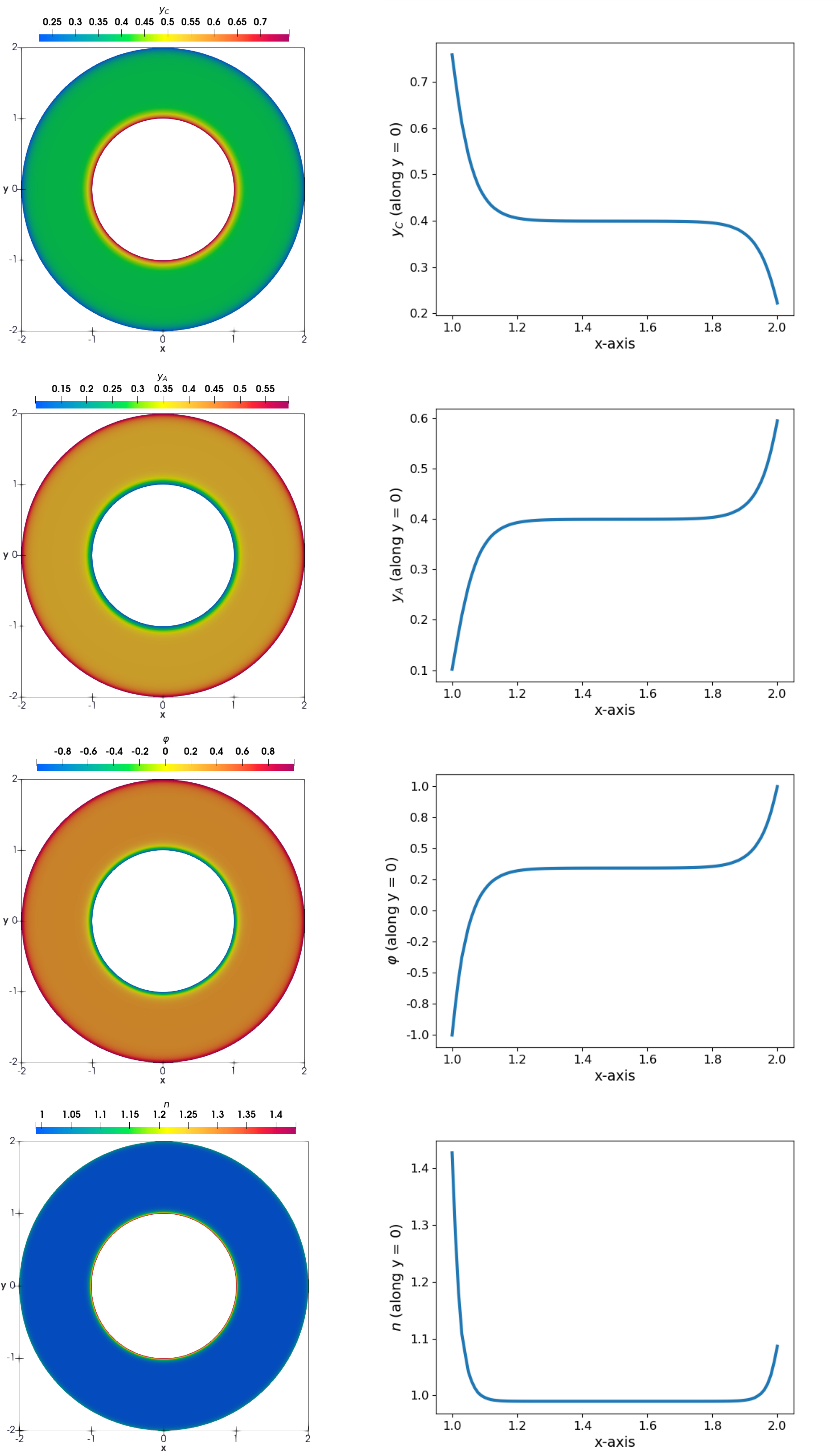}
   \caption{{{Equilibrium distribution of ions (cation and anion), potential, and total number density on full domain (left) and cross-section along $y=0$, $1 \leq x \leq 2$ (right) in an annular battery design in Example~\ref{Sec:6.4.1}, with an applied potential of $-1.0$ unit on the inner boundary and $1.0$ unit on the outer boundary. Parameter settings are provided in Table~\ref{Compressible}. Boundary conditions are described in Section~\ref{Bdry}.}}}
	\label{Fig:8} 
\end{figure}
\begin{figure}[h!]
	\centering
	\includegraphics[width=0.7\linewidth]{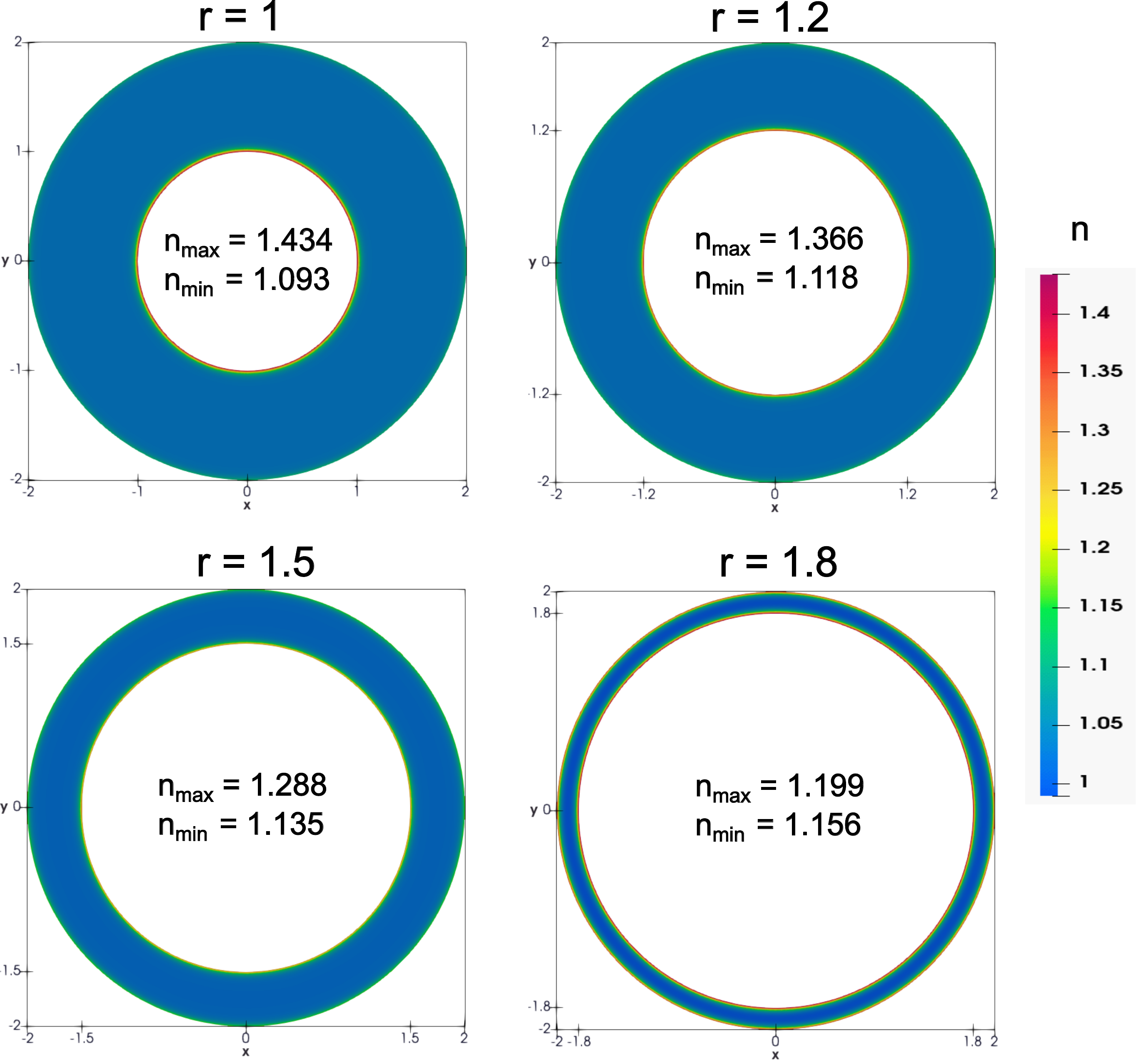}
	\caption{{Numerical solution of total number density with increasing inner radius in example \ref{Sec:6.4.1}. Parameter settings are given in Table~\ref{Compressible}. Boundary conditions are described in Section~\ref{Bdry}, with an applied potential of $-1.0$ unit on the inner boundary and $1.0$ unit on the outer boundary}}  
	\label{Fig:9} 
\end{figure}
\begin{figure}[h!]
	\centering
	\includegraphics[width=0.7\linewidth]{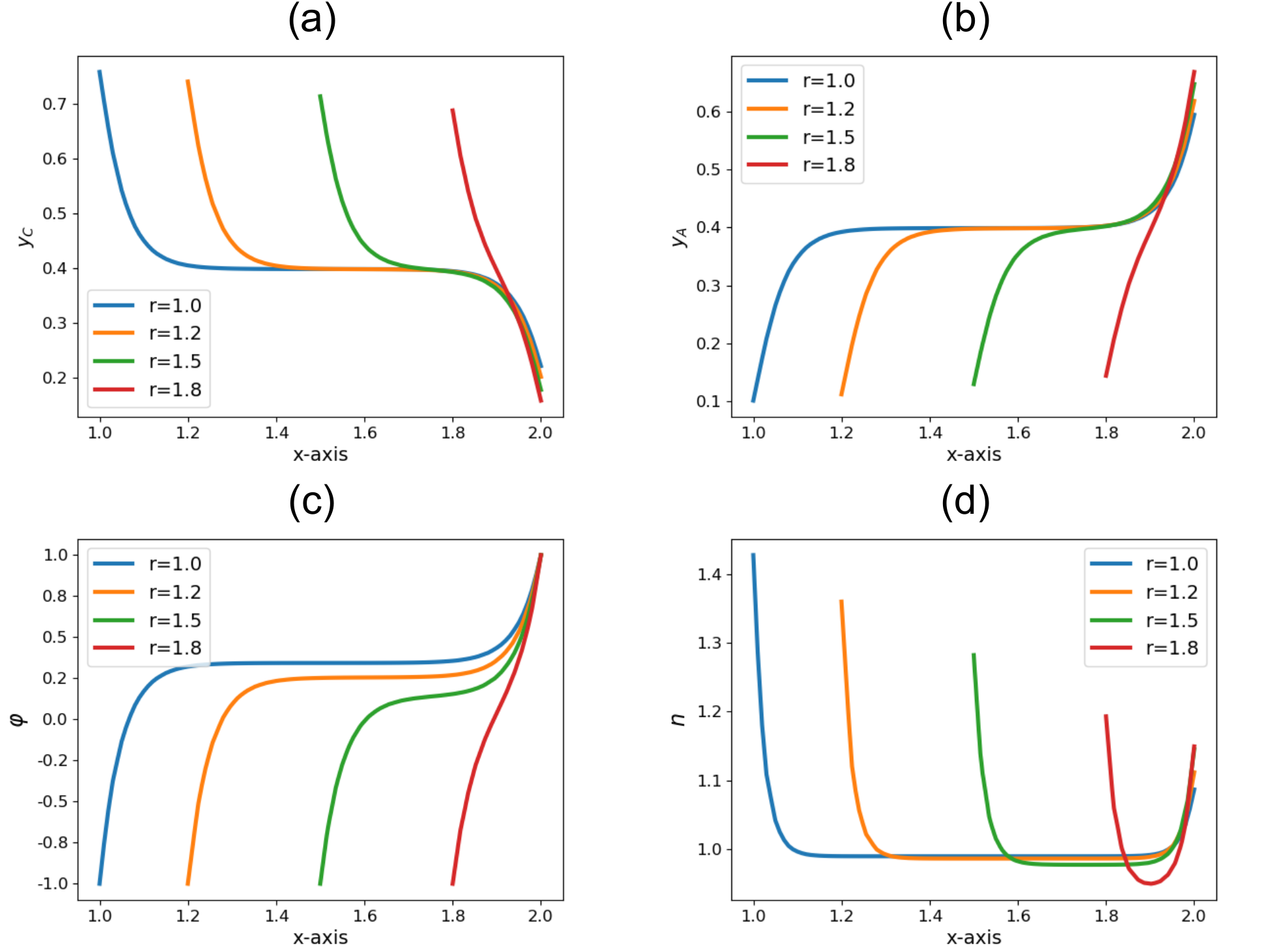}
	\caption{{Numerical solutions of (a) cation $y_c$, (b) anion $y_a$, (c) electric potential $\varphi$, (d) total number density $n$ for annular battery along cross-section $y=0, \;1\leq x\leq 2$ with increasing inner radius in example \ref{Sec:6.4.1}. Parameter settings are given in Table~\ref{Compressible}. Boundary conditions are described in Section~\ref{Bdry}, with an applied potential of  $-1$ unit on the inner boundary and $1$ unit on the outer boundary ($-1$ unit at $x=1$ and $1$ unit at $x=2$).}}
	\label{Fig:10} 
\end{figure}

{Since the system is in equilibrium, we focus on the ionic distributions within the ring-type geometry.} On the inner circle, we prescribe the electrical potential $\varphi=-1.0$ and on the outer circle, we have $\varphi=1.0$. The chosen parameters align with those outlined in Table \ref{Compressible}. Figure \eqref{Fig:8} presents insightful visualizations, revealing a noteworthy asymmetry in the distribution and behavior of cations and anions. This phenomenon can be comprehended by recognizing that the outer circle, being larger, allows for a more extensive dispersion of cations and anions over a larger area, while the inner circle, having less area for dispersion, results in a higher concentration of cations. A parallel non-symmetric behavior is also evident in the variables $\varphi$ and $n$, as illustrated by cross-sectional plotting along $y=0$ in Fig. \eqref{Fig:8}. The observed behavior finds further justification in Fig. \eqref{Fig:9} and \eqref{Fig:10}, where the adjustment of the inner circle's radius contributes to a more symmetrical outcome.  For instance, selecting an inner radius of 1.8, closer to the outer radius of 2, induces a notable shift towards symmetric behavior in the number density. Cross-sectional plotting along $y=0$ in Figure \eqref{Fig:10} provides a comparison of all governing variables for increasing values of the inner radius from 1 to 1.8, leading to the symmetrical solution. It is noteworthy that while the average values remain consistent, the alteration in the inner circle's radius brings about a significant redistribution. For total number density, the higher concentration observed at the anode diminishes, and conversely, the concentration at the cathode increases significantly. 


\subsection{Effect of temperature}
\label{Sec:6.4}
In earlier investigations, we focused on the equilibrium distribution of ions in a mixture under isothermal conditions, as noted by Dreyer et al. \cite{dreyer2013overcoming}.  The current analysis examines the effect of temperature on these equilibrium distributions.  {It should be noted that each simulation is still isothermal, with temperature treated as a parameter.}
It is important to note that temperature varies inversely with the dimensionless parameters
$\Psi, \Lambda$ and $\Hat{K}$.

\begin{figure}[h!]
	\centering
	\includegraphics[width=1.0\linewidth]{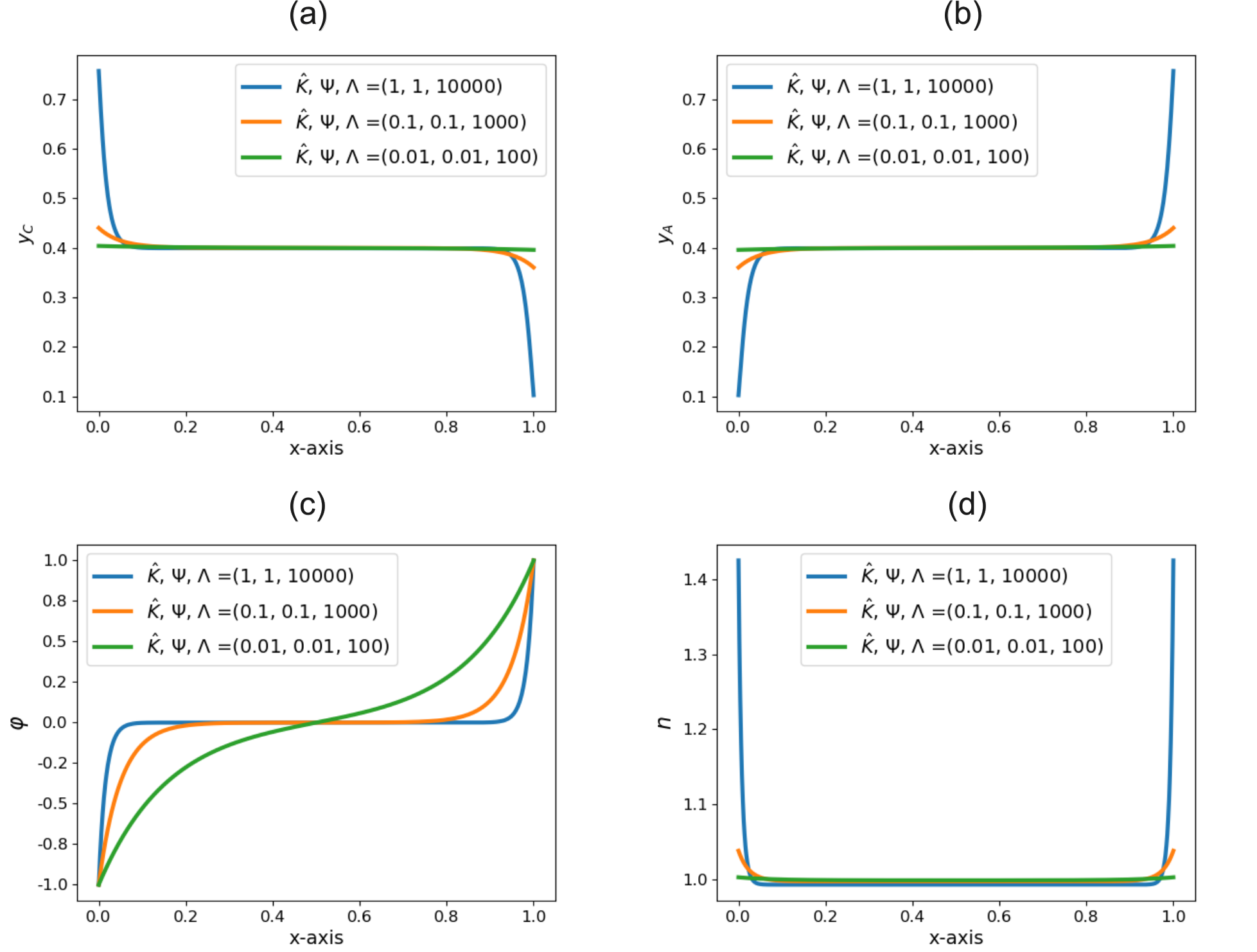}
\caption{{Effect of temperature on numerical solutions in Example~\ref{Sec:6.4}. Parameter values are taken from Table~\ref{Compressible}, except for the varying quantities $\Psi$, $\Lambda$, and $\Hat{K}$. Boundary conditions are described in Section~\ref{Bdry}, with an applied potential of $-1$ unit at $x=0$ and $1$ unit at $x=1$.}}
	\label{Fig:11}  
\end{figure}

The numerical results shown in Fig. \eqref{Fig:11} reveal a notable trend: at lower temperatures, a higher accumulation of cations occurs on the left-hand side, with comparatively fewer on the right-hand side. This observation suggests that, under lower temperature conditions, the influence of the electric potential becomes more pronounced, leading to a preferential attraction of ions towards specific electrodes. A similar influence is observed for other governing variables for low temperature.  

Conversely, as the temperature increases, the diffusive effect becomes more dominant. A higher diffusive effect relative to the electric potential implies that ions are distributed more smoothly across the entire domain rather than being strongly attracted to specific electrodes. This dynamic is evident in Figure \eqref{Fig:11}, where, for elevated temperatures, all the governing variables diffuse throughout the entire domain.
In contrast, for lower temperatures, the attractive forces of the electric potential prevail, resulting in a concentration of ions being drawn towards the electrodes. This stark contrast in ion distribution patterns highlights the significant impact of temperature on the interplay between diffusive and electric potential effects, providing valuable insights into the nuanced behavior of ion transport phenomena under varying thermodynamic conditions. Hence, battery designs should account for these effects to ensure efficient charge transfer.

\begin{figure}[h!]
\centering
 \includegraphics[width=7cm,height=6cm]{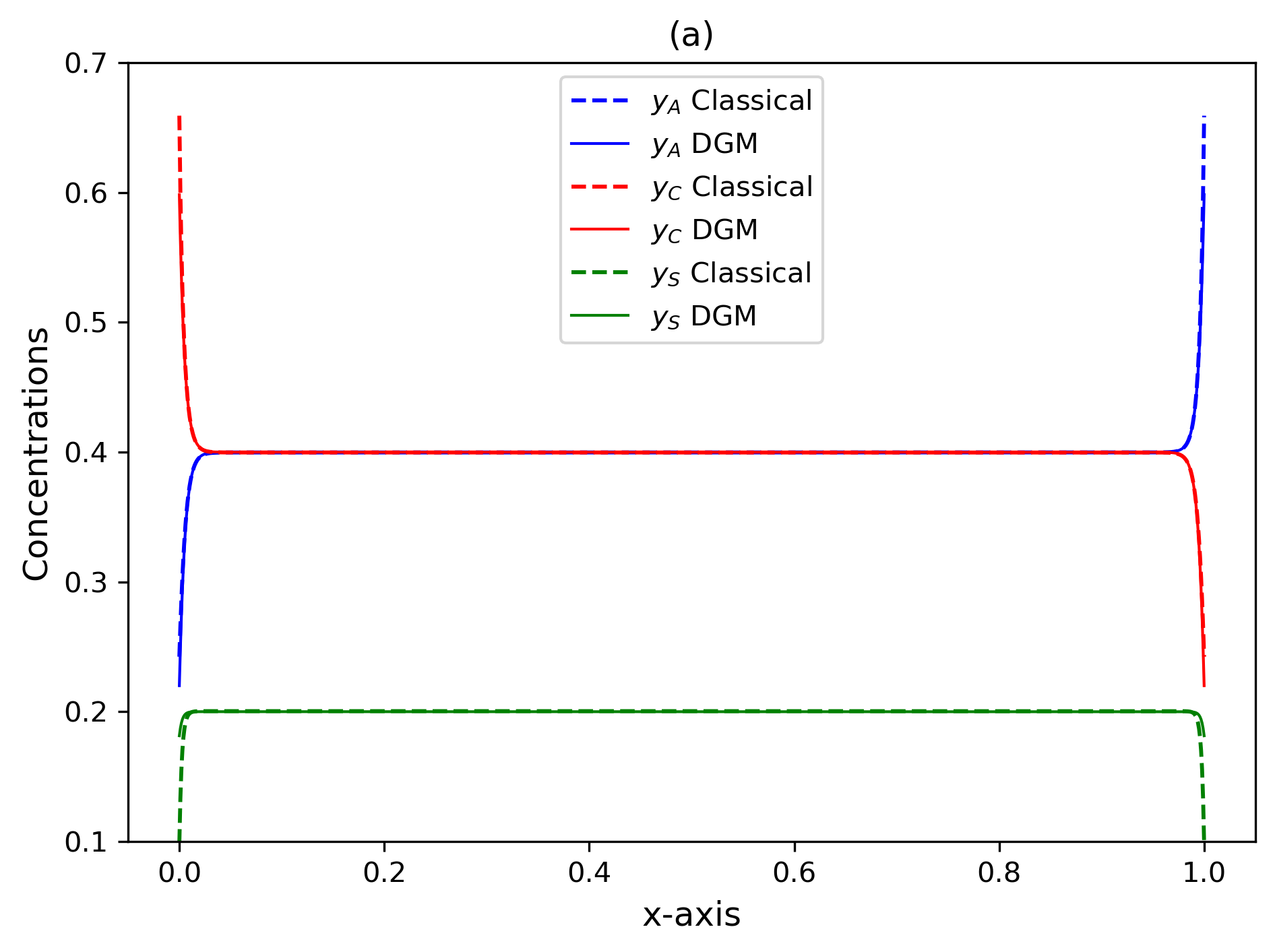}\;\;\hspace{1cm}
  \includegraphics[width=7cm,height=6cm]{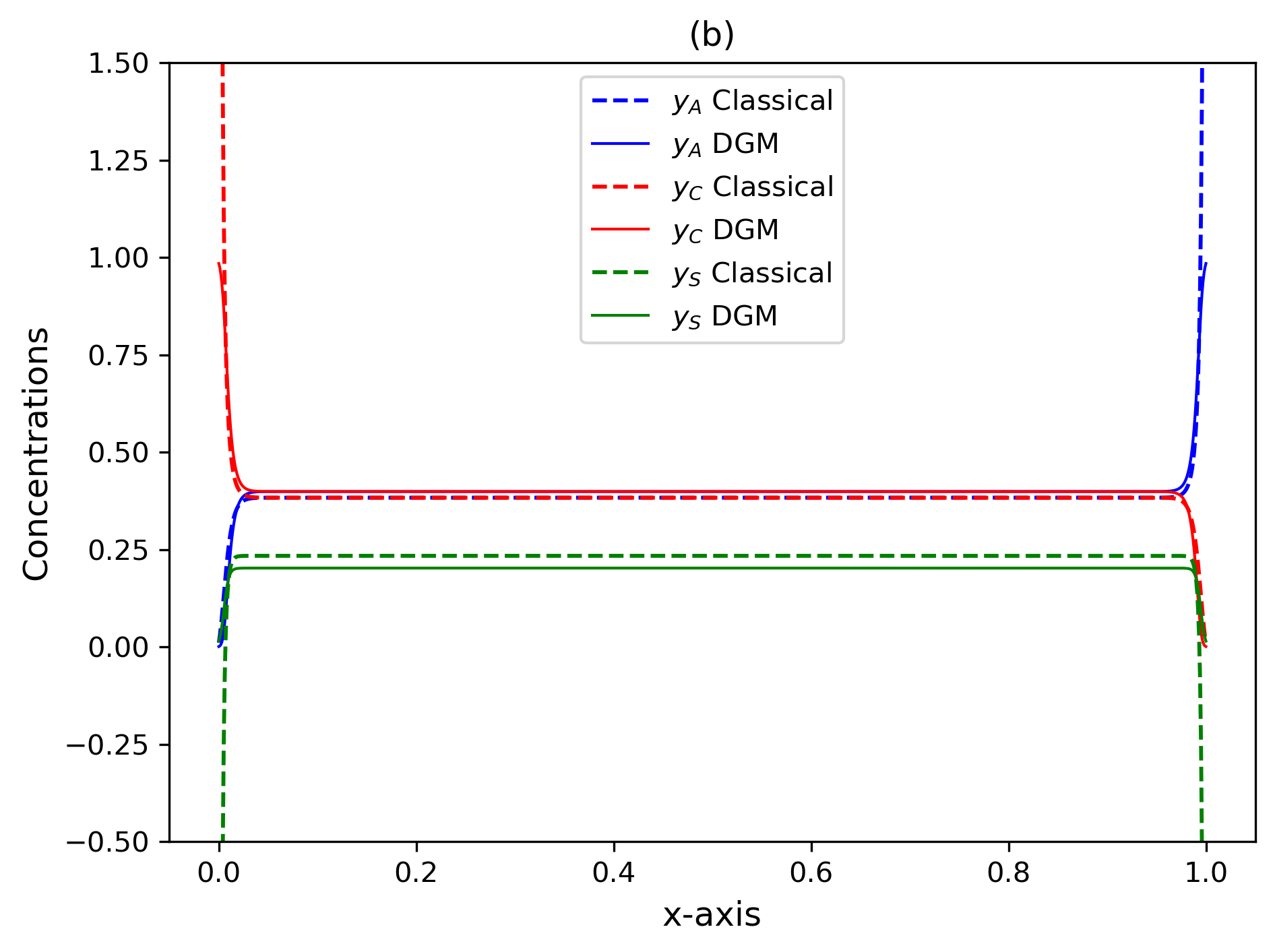}
\caption{{{Comparison of the DGM model with the classical Nernst--Planck model for an incompressible (${K} \rightarrow \infty$), symmetric 1:1 electrolyte (Example~\ref{Sec:6.6}). Simulations use the parameters from Table~\ref{Compressible}, carried over the unit interval with boundary conditions defined in Section~\ref{Bdry}. A dimensionless potential difference of (a) $ 1$ unit, and (b) $ 7$ unit is imposed, with $\phi = 1$ or $7$ at $x = 0$ and $\phi = 0$ at $x = 1$.}}}
\label{Fig:12}  
\end{figure}

\subsection{{Comparison with the classical Nernst-Planck model} } \label{Sec:6.6}

 {This section aims to highlight the key physical differences between the classical PNP model and the thermodynamically consistent Dreyer–Guhlke–Müller (DGM) framework, and to examine how these differences influence numerical behavior, including charge distribution and species transport. Unlike the PNP model, where ionic transport depends solely on each species’ own concentration gradient and the electric field (Eq.~\eqref{Eq8nnn}), the DGM flux (Eq.~\eqref{Eq36}) introduces three key modifications:}

{First, the cross-diffusion terms ($D_{\alpha i}$) account for how gradients in other ionic species ($\nabla y_i$) directly influence the flux of species $\alpha$, capturing important crowding effects in concentrated solutions. This implies that the transport of one species is coupled to the movement of others. Second, the nonlinear self-diffusion coefficient ($D_{\alpha\alpha} \propto 1/y_\alpha$) describes the decreasing mobility of a species as its local concentration increases, reflecting how ions slow down when densely surrounded by others. 
Third, the electromigration term ($D_{\alpha\varphi}$) incorporates local charge density effects through $\sum z_iy_i$, demonstrating how each ion's response to electric fields depends on the net charge distribution of neighboring ions.}

{These modifications become particularly important in extreme regimes. While the classical PNP formulation remains suitable for dilute solutions, the DGM framework provides essential corrections for concentrated electrolytes and other applications where cross-species coupling, strong potential gradients, and nonlinear mobility effects are significant. The DGM model incorporates more realistic ion–ion interactions and accurately captures saturation behavior as concentrations approach their physical bounds (\(0 < y_\alpha < 1\)). In contrast, the classical PNP model may yield nonphysical results under such conditions, which we illustrate below.}

{To begin, we nondimensionalize the classical diffusion flux given in equation~\eqref{Eq8nnn} using the dimensionless parameters defined in Section~\ref{dimensionless}:}
\begin{equation}
\label{Eq8nnnnnn}
{
\bm{J}_{\alpha}= - \frac{M_{\alpha}^{NP} k T n_0}{M_0 x_0} \left(n \nabla y_{\alpha} + y_{\alpha}\nabla n+ \Psi z_{\alpha} y_{\alpha} n \nabla \varphi \right), \quad \alpha \in \{1, 2, \ldots, N\},}
\end{equation}
{here \(M_0\) denotes the reference mobility. We compare the models for an incompressible, symmetric 1:1 electrolyte using the parameters listed in Table~\ref{Compressible}, under both low and high potential differences.  As explained in Section~\ref{Sec:6.3}, in this case, when the bulk modulus \(K \rightarrow \infty\), the total number density \(n\) approaches a constant reference value \(n^R = n^0\). Consequently, the pressure can no longer be determined from the constitutive law \eqref{Eq10n}. Instead, the pressure equation decouples from the system and must be computed separately using equation~\eqref{Eq21*n}, once the remaining governing variables $y_C, y_A$ and $\varphi$ have been determined.}

{The corresponding numerical results are shown in Figure~\ref{Fig:12}. For small potential differences (e.g., a dimensionless value of 1), both models yield nearly identical behavior, with atomic fractions remaining within the physical range \(0 < y_\alpha < 1\). However, as the applied potential difference increases (e.g., to a value of 7), the classical PNP model predicts anion concentrations that exceed the upper bound of one and neutral solvent concentrations that fall below zero, both of which are nonphysical. In contrast, the DGM model respects these physical constraints by driving the atomic fractions toward their saturation limits of zero and one. A similar phenomenon is illustrated in Figure~9 of~\cite{dreyer2013overcoming}, where the classical model produces unrealistically large boundary pressures under strong potential gradients, further highlighting its limitations.}

\subsection{$K$-dependent double layer capacitance curves}
\label{Sec:6.7}
The double layer capacitance is a fundamental property of electrochemical interfaces that quantifies the ability of the interface to store charge in response to an applied potential difference. In this subsection, we analyze how the double layer capacitance varies with the compressibility parameter $K$. We consider a benchmark simulation following \cite{Dreyer2018-1}, extended here to a compressible 1:1 electrolyte with cation charge number $z_A = +1$ and anion charge number $z_C = -1$, at a molarity $M = 0.1~\mathrm{mol/L}$  on the unit interval. The simulation employs the following standard and physically meaningful parameters \cite{dreyer2013overcoming, Dreyer2014-1, Dreyer2018-1}:
\begin{align*}
\varphi_{BC} &= \frac{k T}{e_0}, & 
k &= 1.381 \times 10^{-23} \,\mathrm{J/K}, & 
T &= 293.75 \,\mathrm{K}, & 
e_0 &= 1.602 \times 10^{-19} \,\mathrm{C},\\[2mm]
n^R &= 55 \times N_A \times 10^{-3} \,\mathrm{mol/L}, & 
N_A &= 6.022 \times 10^{23} \,\mathrm{mol^{-1}}, & 
p^R &= 1.01325 \times 10^{5} \,\mathrm{Pa},\\[1mm]
x_0 &= 20 \times 10^{-9} \,\mathrm{m}, & 
\epsilon_0 &= 8.85 \times 10^{-12} \,\mathrm{F/m}.
\end{align*}
In addition to the boundary conditions \eqref{eq51} and \eqref{eq51nn}, we impose Dirichlet boundary conditions at the right boundary $x=1$ as follows:
\begin{equation}
\varphi^{\mathrm{Right}} = 0.0, \qquad p^{\mathrm{Right}} = 0.0, \qquad y_\alpha^{\mathrm{Right}} = \frac{M}{n^{\mathrm{ref}}}, \quad \alpha \in \{A, C\}. 
\end{equation}
With these conditions, the side conditions \eqref{eq50} and \eqref{eq52} are no longer required and we proceed to solve for the fields $y_A$, $y_C$, $\varphi$, and $p$. The {differential double layer capacitance} is defined in terms of the total electrostatic charge stored in the system, which is given by
\begin{equation}
Q(\varphi^{\mathrm{Left}}) = \int_{\Omega} n \left( e_0 z_A y_A + e_0 z_C y_C \right) \, dx.
\end{equation}
The {double layer capacitance} ($C_{dl}$) is then obtained as
\begin{equation}
C_{dl} = \frac{dQ}{d\varphi^{\mathrm{Left}}}.
\end{equation}

\begin{figure}[h!]
\centering
 \includegraphics[width=7cm,height=6cm]{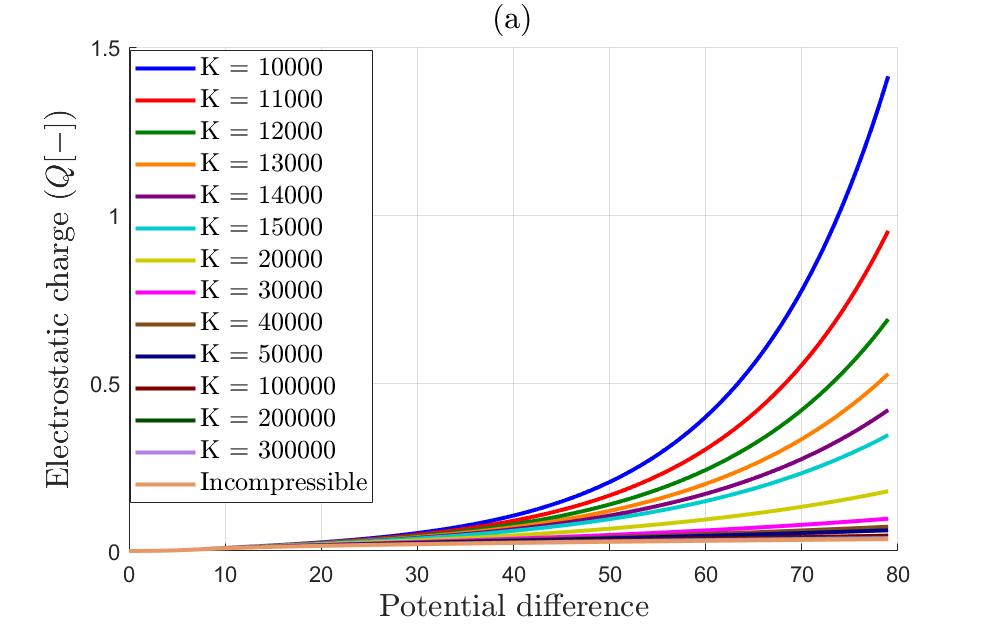}\;\;\hspace{1cm}
  \includegraphics[width=7cm,height=6cm]{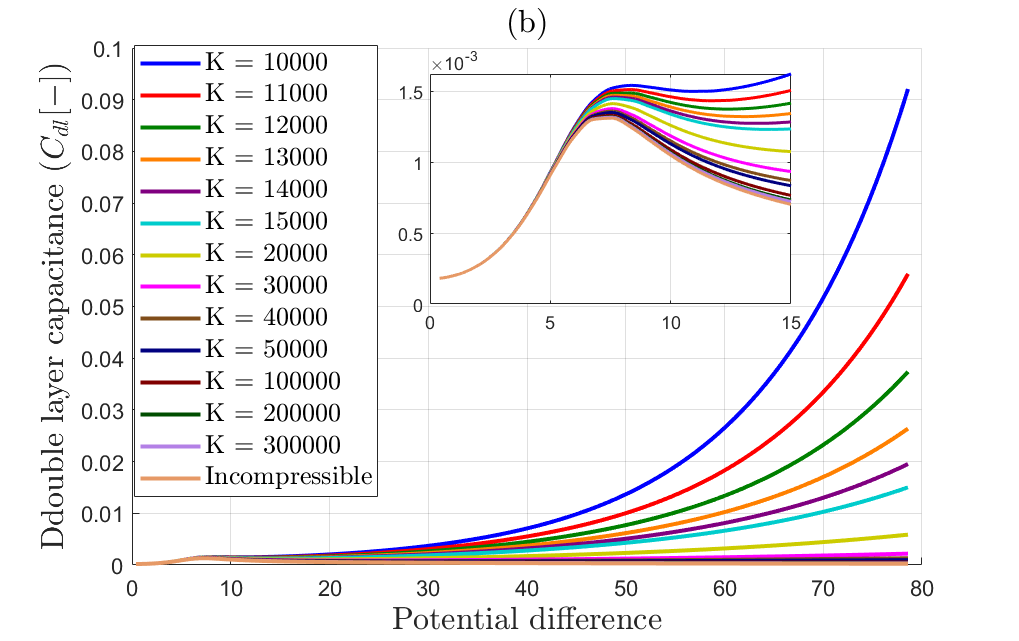}
\caption{{Effect of the compressibility parameter (bulk modulus) $K$ on (a) the total electrostatic charge $Q$ and (b) the differential double layer capacitance $C_{dl}$ for the compressible 1:1 electrolyte system described in Section~\ref{Sec:6.7}.}}
\label{Fig:13}  
\end{figure}

Figure~\ref{Fig:13} illustrates the effect of compressibility on the differential double layer capacitance and the electrostatic charge, each calculated over a wide range of $\varphi^{\mathrm{Left}}$ values. 
The results show that the total electrostatic charge stored in the system increases monotonically as $\varphi^{\mathrm{Left}}$ is varied from 0 to 80 for all values of the compressibility parameter $K$. As $K$ is increased, the total charge decreases and approaches the value corresponding to the incompressible mixture, in agreement with theoretical expectations. This behavior confirms that {decreasing the compressibility of the system} reduces the overall charge accumulation in the double layer. Also, near the electrodes, one ionic species becomes strongly depleted due to the applied potential, resulting in a negligible contribution of that species to the total charge. Consequently, the double layer behavior in this region is primarily governed by the counter-ions that accumulate close to the electrode surface.

The differential double layer capacitance, $C_{dl}$, exhibits a similar dependence on the compressibility parameter $K$, but shows a distinct behavior with respect to the applied potential. 
 For high values of $K$ and the incompressible case, the capacitance exhibits the characteristic camel-shaped profile. For lower values of $K$, the capacitance initially increases with $\varphi^{\mathrm{Left}}$, reaches a shallow local maximum, then decreases slightly before increasing again at higher potentials.
 This non-monotonic behavior is associated with changes in the charge distribution and may be influenced by the additional density variations permitted in the compressible electrolyte (see Section~\ref{Sec:6.3}).
  As $K$ increases, the overall magnitude of $C_{dl}$ decreases and approaches the incompressible limit, indicating that {less compressible systems (high values of $K$)} have a reduced ability to store charge within the double layer. 
  These results highlight that electrolyte compressibility can play an important role in regimes where density variations significantly affect interfacial charge distributions.

\section{Concluding remarks and outlook}
\label{Sec:6}

This study offered a finite element framework for the multidimensional simulation of modified Poisson–Nernst
–Planck/Navier–Stokes (PNP/NS) system. The modification involves substituting the diffusion flux in the original Nernst–Planck equation with an implicitly defined diffusion flux \cite{dreyer2013overcoming}. To develop the comprehensive framework, we utilized the features offered by the open-source software FEniCS.
Through rigorous numerical validation, we demonstrated the accuracy and reliability of the proposed approach. Our investigation encompassed a meticulous examination of both compressible and incompressible mixtures, shedding light on parametric dependencies governing space-charge layer formation. Moreover, the present study extends to the nuanced behavior of ions under diverse operational conditions, including temperature fluctuations and the influence of bulk modulus.  

Looking ahead, there is potential for extending our structured framework to encompass non-equilibrium mixtures. Furthermore, the forthcoming work will delve into a detailed analysis, including the demonstration of well-posedness and derivation of error estimates in both semi-discrete and fully discrete formulations.

\section*{Declaration of competing interest}
 The authors declare that they have no conflict of interest.

\section*{Acknowledgments}
The authors are grateful to the anonymous referees for their valuable comments and suggestions

\vspace{0.5cm}

Ankur acknowledges the funding through the Advanced Research Opportunities Program (AROP), Aachen University, Germany. 
S.S. is funded by German Research Foundation (DFG) under Research Unit FOR5409: "Structure-Preserving Numerical Methods for Bulk- and Interface-Coupling of Heterogeneous Models (SNuBIC)".
The authors also sincerely thank Prof. Manuel Torrilhon (Head of the ACoM Lab, RWTH Aachen University, Germany) for his valuable guidance and intellectual input throughout the entire AROP project carried out by Ankur at the ACoM Lab.

\section*{Data availability}
{The implementation code and simulation data supporting the findings of this study are available at the GitHub repository: \url{https://github.com/Ankur-IIT/Modified_PNP-NS_Model}}.

%
%
%
%
%
%

\appendix
\section{Mathematical formulation for modified PNP/NS model}
\label{Sec:2app}
Let us consider an ideal elastic mixture composed of $N$ constituents $I_\alpha$, $\alpha \in \{1, 2, \ldots N\}$. Among these constituents, $I_1, I_2,$ $\ldots I_{N-1}$  are ionic species and $I_N$ is neutral solvent.
The evolution of this system is governed by $N-1$ partial mass balances, the Poisson equation and Navier–Stokes equations incorporating a body force controlled by the electric field and the free charge density. To account for the finite ion size effects in the constitutive relationship, the authors of \cite{dreyer2013overcoming,de2013non} suggested that chemical potentials having a relationship with concentrations, atomic masses and pressure  
have to be considered as a driving force. Two major applications for these kinds of frameworks are in semiconductor devices and electrochemistry. 

\subsection{Mass conservation equations} 
\label{Sec:2.1app}
Suppose the mixture occupies a spatial domain $\Omega\in \mathbb{R}^n, n\leq 3$ such that there is no flux through the boundary. Also, assume that the constituents of the mixture have atomic masses ($m_{\alpha})_{\alpha \in \{1, 2, \ldots N\}}$, charge numbers ($z_{\alpha})_{\alpha \in \{1, 2, \ldots N\}}$,  number densities ($n_{\alpha})_{\alpha \in \{1, 2, \ldots N\}}$ and velocities ($\bm{v}_{\alpha})_{\alpha \in \{1, 2, \ldots N\}}$ at any time $t\geq 0$. Using these notations, we define the partial mass densities 
\begin{equation}
\label{Eq1}
\rho_{\alpha}= m_{\alpha}n_{\alpha}, \hspace{1cm} \alpha \in \{1, 2, \ldots N\},
\end{equation}
total mass density  
\begin{equation}
\label{Eq2}
\rho= \sum_{\alpha=1}^{N}\rho_{\alpha},
\end{equation}
total number density  
\begin{equation}
\label{Eq2new}
n= \sum_{\alpha=1}^{N}n_{\alpha},
\end{equation}
and the barycentric velocity as
\begin{equation}
\label{Eq3}
\bm{v}= \frac{1}{\rho}\sum_{\alpha=1}^{N}\rho_{\alpha}\bm{v}_{\alpha}.
\end{equation}
The diffusive mass flux in reference to the constituents $I_{\alpha}$  is given by 
\begin{equation}
\label{Eq4}
\begin{aligned}
\mbox{Diffusion\; flux} &= \mbox{Total\; flux} - \mbox{convective\; flux} \\
\bm{J}_{\alpha}~~~~&=~~~~~\rho_{\alpha}\bm{v}_{\alpha}~~~~-~~~\rho_{\alpha}\bm{v}. 
\end{aligned}
\end{equation}
One may note that $\sum_{\alpha=1}^{N}\bm{J}_{\alpha}=0.$
Utilizing the partial mass balance equation for each constituent of the mixture, characterized by densities $\rho_{\alpha}$, we obtain the fundamental relation
\begin{equation}
\label{Eq5}
\frac{\partial \rho_{\alpha}}{\partial t} + \nabla\cdot(\rho_{\alpha}\bm{v}_{\alpha})=0.
\end{equation}
Further, incorporating equation \eqref{Eq4}, we arrive at
\begin{equation}
\label{Eq6}
\frac{\partial \rho_{\alpha}}{\partial t} + \nabla\cdot(\rho_{\alpha}\bm{v}+\bm{J}_{\alpha})=0.
\end{equation}
By summing up the above equation for all constituents i.e. $\alpha=1,2\ldots N$ and leveraging equation \eqref{Eq2} along with the crucial condition $\sum_{\alpha=1}^{N}\bm{J}_{\alpha}=0$, we get the total mass balance as
\begin{equation}
\label{Eq7}
\frac{\partial \rho}{\partial t} + \nabla\cdot(\rho\bm{v})=0.
\end{equation}
The thermodynamic consistent PNP/NS model is based on these $N-1$ partial mass balances (equation \ref{Eq6}) of ionic species and the total mass balance equation \eqref{Eq7}. In conjunction with these partial mass balances, the model is further characterized by two additional equations, which we will derive below.
\subsection{Momentum conservation equation and Poisson equation} \label{Sec:2.2app} 
Utilizing the notation introduced in \ref{Sec:2.1app}, the free charge density of the mixture is expressed as 
\begin{equation}
\label{Eq11}
n^F= \sum_{\alpha=1}^{N}z_{\alpha}e_0n_{\alpha}, \hspace{1cm} \alpha \in \{1, 2, \ldots N\},
\end{equation}
where $e_0$ represents the elementary charge. Alongside the free charge density, additional charge densities arise from polarization effects. The total charge density ($n^e$) will be the sum of free charge density and charge density due ($n^P$) to the polarization. Using Maxwell's theory \cite{Muller1985}, the expression for $n^P$ in terms of a polarization vector $\bm{P}$ is given by
\begin{equation}
\label{Eq46}
n^P= -\nabla\cdot\bm{P}.
\end{equation}
The conservation equation for momentum is expressed as
\begin{equation}
\label{Eq12}
\frac{\partial (\rho\bm{v})}{\partial t} + \nabla\cdot(\rho\bm{v}	\otimes\bm{v}-\bm{\sigma})=\rho\bm{b}+\bm{k}.
\end{equation}
Here, $\bm{\sigma}$ denotes the stress tensor, and $\rho\bm{b}$ and $\bm{k}$ represent the force densities due to gravitation and electromagnetic fields. We consider only the contribution from the electric field ($\bm{E}$), leading to
\begin{align}
\bm{b}&=  0,\label{Eq13}\\
\bm{k}&=n^e\bm{E}.\label{Eq14} 
\end{align}
Using Maxwell's equations for quasi-static electric fields, we have 
\begin{equation}
\label{Eq15}
\bm{E} = - \nabla \varphi, \quad n^e=\epsilon_0 \nabla\cdot(\bm{E}).
\end{equation}
Here, $\epsilon_0$ and $\varphi$ denote the dielectric constant and electric potential respectively. Thus, Lorentz force density ($\bm{k}$) becomes
\begin{equation}
\label{Eq16}
\bm{k}=\epsilon_0 \nabla\cdot(\bm{E})\bm{E}=\nabla\cdot\bigg(\epsilon_0 (\bm{E}\otimes\bm{E})-\frac{1}{2}\epsilon_0\left|\bm{E}\right|^2\bm{1}\bigg).
\end{equation}
Additionally, for the given dielectric susceptibility $\chi$, the stress tensor $\bm{\sigma}$ and polarization vector $\bm{P}$ are given by \cite{dreyer2013overcoming}
\begin{subequations}
	\label{Eq17}
	\begin{align}
	\bm{\sigma}&=- p\bm{1}-\frac{1}{2}\epsilon_0\chi\left|\bm{E}\right|^2\bm{1}+ \epsilon_0\chi (\bm{E}\otimes\bm{E}),\label{Eq47}\\
	\bm{P}&=\epsilon_0\chi\bm{E}, \label{Eq48}
	\end{align}
\end{subequations} 
where $p$ and $\bm{1}$ represent the pressure and unit matrix respectively. Finally, using equations \eqref{Eq46}, \eqref{Eq13},$\eqref{Eq15}$, \eqref{Eq16} and \eqref{Eq17}, the system is reduced to following coupled system of momentum balance and Poisson's equation 
\begin{align}
\frac{\partial (\rho\bm{v})}{\partial t} + \nabla\cdot(\rho\bm{v}	\otimes\bm{v})+ \nabla p&=-n^F\nabla\varphi, \label{Eq18a}\\
-\epsilon_0\epsilon_r\Delta\varphi &= n^F, \label{Eq18b}
\end{align}
where $\epsilon_r=1+\chi$, represents relative dielectric permittivity. Thus, equations \eqref{Eq6} for $\alpha=1,2, \ldots N-1$ together with \eqref{Eq7}, \eqref{Eq18a} and \eqref{Eq18b} represent our system.  The full set of equations for the modified PNP/NS model reads as 
\begin{subequations}
	\label{Eq23}
	\begin{align}
	&\frac{\partial \rho_{\alpha}}{\partial t} + \nabla\cdot(\rho_{\alpha}\bm{v}+\bm{J}_{\alpha})=0, \quad \alpha\;\in\;\{1,\ldots,N-1\}, \label{Eq19}\\
	&\frac{\partial \rho}{\partial t} + \nabla\cdot(\rho\bm{v})=0,\label{Eq20}\\
	&\frac{\partial (\rho\bm{v})}{\partial t} + \nabla\cdot(\rho\bm{v}	\otimes\bm{v})+ \nabla p=-n^F\nabla\varphi,\label{Eq21}\\
	&-\epsilon_0\epsilon_r\Delta\varphi = n^F. \label{Eq22}
	\end{align}
\end{subequations} 

\subsection{Closure of the system} \label{Sec:2.3app}
Up to now, we have compiled a comprehensive set of equations to model the behavior of the  mixture. In accordance with the findings presented in \cite{dreyer2013overcoming}, for a given temperature $T$ and elementary charge $e_0$ , the diffusion fluxes in equation \eqref{Eq19} can be expressed as follows  
\begin{equation}
\label{Eq8}
\bm{J}_{\alpha}= - \sum_{\beta=1}^{N-1}H_{\alpha\beta} \bigg(\frac{\nabla\mu_{\beta}-\nabla\mu_N}{T}+ \frac{e_0}{T}\bigg(\frac{z_{\beta}}{m_{\beta}}-\frac{z_{N}}{m_{N}}\bigg)\nabla\varphi\bigg), \quad\alpha \in \{1, 2, \ldots N-1\}, 
\end{equation}
where $\varphi$, $\mu_{\alpha}$'s and $H_{\alpha\beta}$  denote the electric potential, chemical potentials, and positive definite kinetic matrix, respectively. 

As we have considered the ideal elastic mixture, the chemical potentials are intricately related to concentrations, atomic masses, and specific Gibbs energies $(g_{\alpha})$ through the following expression 
\begin{equation}
\label{Eq9}
\mu_{\alpha}= g_{\alpha}^R+\frac{kT}{m_{\alpha}}\mbox{ln}\bigg(\frac{n_\alpha}{n}\bigg)+\frac{K}{m_{\alpha}n^R}\mbox{ln}\bigg(1+\frac{p-p^R}{K}\bigg) , \quad\quad\alpha \in \{1, 2, \ldots N\}. 
\end{equation}
Further, material pressure is related to total number density as follows
\begin{equation}
\label{Eq10}
p= p^R+K\bigg(\frac{n}{n^R}-1\bigg), 
\end{equation}
where reference state is denoted by index $R$ and  Boltzmann constant by $k$. We use the bulk modulus $K$ to introduce the notion of compressibility in the later half of the article.

\section{Notations}
\label{Notations}
\begin{longtable}{|c|c|c|}
	\captionsetup{font=footnotesize}
	\caption{ List of Symbols} \label{tab1:symbols} \\
	\hline
	\textbf{Symbol} & \textbf{Definition} & \textbf{Relation with other quantities} \\
	\hline
	\endfirsthead
	\multicolumn{3}{c}%
	{\tablename\ \thetable{} -- Continued from previous page} \\
	\hline
	\textbf{Symbol} & \textbf{Definition} & \textbf{Relation with other quantities} \\
	\hline
	\endhead
	\hline \multicolumn{3}{r}{Continued on next page} \\
	\endfoot
	\hline
	\endlastfoot
	
	$m_{\alpha}$ & Atomic mass of constituent $I_{\alpha}$ & $-$ \\
	$z_{\alpha}$ & Charge number of constituent $I_{\alpha}$ & $-$ \\
	$n_{\alpha}$ & Number density of constituent $I_{\alpha}$ & $-$ \\
	$\bm{v}_{\alpha}$ & Velocity of constituent $I_{\alpha}$ & $-$ \\
	$\rho_{\alpha}$ & Partial mass density of constituent $I_{\alpha}$ & $\rho_{\alpha}= m_{\alpha}n_{\alpha}$ \\
	$\rho$ & Total mass density of mixture & $\rho= \sum_{\alpha=1}^{N}\rho_{\alpha}$ \\
	$n$ & Total number density of mixture & $n= \sum_{\alpha=1}^{N}n_{\alpha}$ \\
	$\bm{v}$ & Barycentric velocity of mixture & $\bm{v}=\frac{1}{\rho}\sum_{\alpha=1}^{N}\rho_{\alpha}\bm{v}_{\alpha}$ \\
	$\bm{J}_{\alpha}$ & Diffusion flux of constituent $I_{\alpha}$ & $-$ \\
	$H_{\alpha\beta}$ & Kinetic matrix  & $-$ \\
	$\varphi$ & Electric potential & $-$ \\
	$e_0$ & Elementary charge  & $-$ \\
	$T$ & Temperature  & $-$ \\
	$K$ & Bulk modulus  & $-$ \\
	$k$ & Boltzmann constant   & $-$ \\
	$*^R$ & Reference state of $*$   & $-$ \\
	$g_{\alpha}$ & Specific Gibbs energy of $I_{\alpha}$& $g_{\alpha}= g_{\alpha}^R+\frac{K}{m_{\alpha}n^R}\mbox{ln}\big(\frac{n}{n^R}\big)$ \\
	$\mu_{\alpha}$ & Chemical potential of $I_{\alpha}$& $\mu_{\alpha}= g_{\alpha}+ \frac{kT}{m_{\alpha}}\mbox{ln}\big(\frac{n_{\alpha}}{n}\big)$ \\
	$p$ & Pressure   & $p=p^R+K\big(\frac{n}{n^R}-1\big)$ \\
	$n^F$ & Free charge density &  $n^F= \sum_{\alpha=1}^{N}z_{\alpha}e_0n_{\alpha}$\\
	$\epsilon_0$ & Dielectric constant &  $-$\\
	$\chi$ & Dielectric susceptibility &  $-$\\
	$\bm{P}$ & Polarization vector & $\bm{P}=\epsilon_0\chi\bm{E}$ \\
	$n^P$ & Charge density due to polarization &  $n^P=-\nabla \cdot\bm{P}$\\
	$\bm{E}$ & Electric field &$-\nabla \varphi$ \\
	$n^e$ & Total charge density &  $n^e= n^F+n^P$\\
	$\epsilon_r$  & Relative dielectric permittivity &  $\epsilon_r=1+\chi$\\
	$\bm{\sigma}$ & Stress tensor &  $\bm{\sigma}=- p\bm{1}-\frac{1}{2}\epsilon_0\chi\left|\bm{E}\right|^2\bm{1}+ \epsilon_0\chi \bm{E}\otimes\bm{E}$\\
	$\bm{k}$  &  Lorentz force density &  $\bm{k}=n^e\bm{E}$\\
	$y_{\alpha}$  &  Atomic fractions &  $y_{\alpha}=\frac{n_\alpha}{n}$\\
	${M}_\alpha$  &  Mass fractions &  ${M}_\alpha= \frac{m_\alpha}{m_N}$\\
\end{longtable}




\end{document}